\documentclass{aoscm}
\usepackage[dvips]{graphicx}
\usepackage{epsfig}
\def\theequation{\arabic{equation}}
\def\eps{\varepsilon}
\def\bbox{\quad\hbox{\vrule \vbox{\hrule \vskip2pt \hbox{\hskip2pt
\vbox{\hsize=1pt}\hskip2pt} \vskip2pt\hrule}\vrule}}
\def\lessim{\ \lower4pt\hbox{$
\buildrel{\displaystyle <}\over\sim$}\ }
\def\gessim{\ \lower4pt\hbox{$\buildrel{\displaystyle >}
\over\sim$}\ }

\def\eps{{\varepsilon}}

\def\goin{\to\infty}
\newtheorem{proposition}{Proposition}
\newtheorem{lemma}{Lemma}
\newtheorem{theorem}{Theorem}

\makeatletter
\@addtoreset{equation}{section}
\renewcommand{\theequation}{\thesection.\arabic{equation}}
\makeatother

\font\tencmmib=cmmib10 \skewchar\tencmmib '60
\newfam\cmmibfam
\textfont\cmmibfam=\tencmmib

\font\tenmsb=msbm10
\font\eightmsb=msbm10 scaled 800
\def\Bbb#1{\hbox{\tenmsb#1}}
\def\Bbbb#1{\hbox{\eightmsb#1}}
\def\bbox{\quad\hbox{\vrule \vbox{\hrule \vskip2pt \hbox{\hskip2pt
\vbox{\hsize=1pt}\hskip2pt} \vskip2pt\hrule}\vrule}}
\def\lessim{\ \lower4pt\hbox{$
\buildrel{\displaystyle <}\over\sim$}\ }
\def\gessim{\ \lower4pt\hbox{$\buildrel{\displaystyle >}
\over\sim$}\ }

\def\eps{\varepsilon}

\def\goin{\to \infty}
\def\go0{\to 0}

\def\leftitem#1{\item{\hbox to\parindent{\enspace#1\hfill}}}

\def\mdsk{\medskip}

\def\qed{\hfill\break\rightline{$\bbox$}}

\def\sg{\sigma}

\def\sg2{\sigma^2}

\def\__{_{\infty}}

\begin{document}

\title {EMPIRICAL MARGIN DISTRIBUTIONS AND BOUNDING THE GENERALIZATION ERROR
OF COMBINED CLASSIFIERS}
\author{
V. Koltchinskii\thanks{Partially supported by NSA Grant
MDA904-99-1-0031}
\ and D. Panchenko\thanks{Partially supported by Boeing Computer Services Grant 3-48181}
\\
Department of Mathematics and Statistics
\\
The University of New Mexico
}
\maketitle



\vskip 1mm

\centerline{Dedicated to A.V. Skorohod on his seventieth birthday}

\vskip 1mm

\begin{abstract}
We prove new probabilistic upper bounds on generalization error of complex classifiers that
are combinations of simple classifiers. Such combinations could be implemented by neural networks or
by voting methods of combining the classifiers, such as boosting and bagging. The bounds are
in terms of the empirical distribution of the margin of the combined classifier. They are
based on the methods of the theory of Gaussian and empirical processes (comparison inequalities,
symmetrization method, concentration inequalities) and they improve previous results
of Bartlett (1998)
on bounding the generalization error of neural networks in terms of $\ell_1$-norms
of the weights of neurons
and of Schapire, Freund, Bartlett and Lee (1998) on bounding the
generalization error of boosting. We also obtain rates of convergence
in L\'evy distance of
empirical margin distribution to the true margin distribution uniformly over the
classes of classifiers and prove the optimality of these rates.
\end{abstract}

\vskip 5mm

\hfill\break
\hfill\break
\hfill\break
{\it Short title}: Empirical Margins and Generalization Error

\vfill\break

\section{Introduction}

Let $(X,Y)$ be a random couple, where $X$ is an instance in a space $S$ and $Y\in \{-1,1\}$ is a
label. Let ${\cal G}$ be a set of functions from $S$ into ${\Bbb R}.$ For $g\in {\cal G},$
${\rm sign}(g(X))$ will be used as a predictor (a classifier) of the unknown label $Y.$ If the
distribution
of $(X,Y)$ is unknown, then the choice of the predictor is based on the training data
$(X_1,Y_1),\dots ,(X_n,Y_n)$ that consists of $n$ i.i.d. copies of $(X,Y).$ The goal of
learning is to find a predictor $\hat g\in {\cal G}$ (based on the training data) whose
\it generalization (classification) error \rm ${\Bbb P}\{Y\hat g(X)\leq 0\}$ is small enough.
In this paper, our main concern is to find reasonably good probabilistic upper
bounds on the generalization error. The standard approach to this problem was developed
in seminal papers of Vapnik and Chervonenkis in the 70s and 80s (see Vapnik (1998),
Devroye, Gy\" orfi and Lugosi (1996), Vidyasagar (1997)) and it is based on
bounding the difference between the generalization error
${\Bbb P}\{Yg(X)\leq 0\}$ and the training error
$$
n^{-1}\sum_{j=1}^n I_{\{Y_j g(X_j)\leq 0\}}
$$
uniformly over the whole class ${\cal G}$ of classifiers $g$. These bounds are expressed
in terms of data dependent entropy characteristics of the class of sets
$\{\{(x,y): yg(x)\leq 0\}: g\in {\cal G}\}$ or, frequently, in terms
of the so called VC-dimension of the class. It happened, however, that in many important
examples (for instance, in neural network learning) the VC-dimension of the class can be very large,
or even infinite, and that makes impossible the direct application of Vapnik--Chervonenkis
type of bounds.
Recently, several authors (see Bartlett (1998),
Schapire, Freund, Bartlett and Lee (1998), Anthony and Bartlett (1999))
suggested another class of upper bounds on generalization error that are
expressed in terms of the empirical distribution of \it the margin \rm of the predictor
(the classifier).
The margin is defined as the product $Y\hat g(X).$ The bounds
in question are especially useful in the case of the classifiers that are the
combinations of simpler classifiers (that belong, say, to a class ${\cal H}$). One of the examples of
such classifiers is provided by neural networks.
Other examples are given by the classifiers obtained by boosting, bagging
and other voting methods of combining the classifiers.
The bounds in terms of margins are also of interest in application
to generalization performance of support vector machines,
Cortes and Vapnik (1995), Vapnik (1998), Bartlett and
Shawe-Taylor (1999).
The upper bounds have the following form (up to some extra terms)
$$
\inf_{\delta >0} \Bigl[n^{-1}\sum_{j=1}^n I_{\{Y_j\hat g(X_j)\leq \delta\}}+ C({\cal G})
\phi(\delta)\frac{C({\cal H})}{\sqrt{n}}\Bigr],
$$
where $C({\cal G})$ is a constant depending on the class ${\cal G}$
(in other words, on the method of combining the simple classifiers),
$\phi $ is a decreasing function such that $\phi (\delta)\goin$ as $\delta \to 0$
(often, for instance, $\phi (\delta)=\frac{1}{\delta}$), $C({\cal H})$
is a constant depending on the class ${\cal H}$ (in particular, on the VC-dimension,
or some type of entropy characteristics of the class).

It was observed in experiments that classifiers produced by such
methods as boosting tend to have rather \it large margin \rm
of correctly classified examples. This allows one to choose
a relatively large value of $\delta$ in the above bound
without increasing substantially the value of the empirical
distribution function of the margin (which is the first term
of the bound) comparing with the training error.
For large enough $\delta ,$ the second term
becomes small, which ensures a reasonably small value of
the infimum. This allowed the above mentioned authors
to explain partially (at least at qualitative level)
a very good generalization performance of voting and some
other methods of combining simple classifiers observed in many
experiments. This also motivated the development of the methods
of combining the classifiers based on explicit optimization
of the penalized average cost function of the margins, see
Mason, Bartlett and Baxter (1999), Mason, Baxter, Bartlett
and Frean (1999).

Despite the fact that previously developed bounds provide
some explanations of the generalization performance of
complex classifiers, it was actually acknowledged by
Bartlett (1998), Schapire, Freund, Bartlett and Lee (1998)
that the bounds in question have not reached their final
form yet and more research is needed to understand better
the probabilistic nature of these bounds. This becomes
especially important because of the growing number
of boosting type methods (see Friedman, Hastie, Tibshirani (2000),
Friedman (1999)) for which a comprehensive theory is yet to be
developed.
The methods of proof developed by Bartlett (1998) are based on
the so called fat-shattering dimensions of function classes
and on the extension of Vapnik--Chervonenkis type inequalities to such
dimensions. The method of Schapire, Freund, Bartlett and Lee (1998)
exploits the fact that the complex classifiers are \it convex
combinations \rm of base classifiers (these authors
suggest also an extension of their method to the classes
of functions for which there exist so called $\eps$-sloppy
$\theta$-covering). The use of these methods in the case
of general cost functions of the margins poses some difficulties
(see Mason, Bartlett and Baxter (1999)).

In this paper, we develop a new approach that allows us to improve
and better understand some of the previously known bounds.
Our method is based on the general results of the theory
of Gaussian, Rademacher and empirical processes (such as comparison inequalities,
e.g. Slepian's Lemma, symmetrization and random multipliers inequalities,
concentration inequalities, see Ledoux and Talagrand (1991),
van der Vaart and Wellner (1996), Dudley (1999)).
We give the bounds in terms of general functions of the
margins, satisfying a Lipschitz condition.
They can be readily applied to the classifiers based on
explicit optimization of margin cost functions
(such as in the paper of Mason, Bartlett and Baxter (1999)).
In the case of Bartlett's
bounds for feedforward neural networks in terms of the $\ell_1$-norms of the
weights of the neurons (see Bartlett (1998) and also Fine (1999)),
the improvement we got is substantial. In Bartlett's bounds
the constant $C({\cal G})$ is of the order $(AL)^{l(l+1)/2},$
where $A$ is an upper bound on the $\ell_1$-norms of the weights of neurons,
$L$ is the Lipschitz constant of the sigmoids, and $l$ is the number of layers
of the network. Also, in his bound $\phi (\delta)=\frac{1}{{\delta}^l}.$
We obtained in a similar context $C({\cal G})$ of the order $({AL})^l$ with
$\phi (\delta)=\frac{1}{\delta}.$

Based on our bounds, we developed a method
of complexity penalization of the training error of neural network learning with penalties defined
as functionals of the weights of neurons and prove oracle inequalities
showing some form of optimality of this method.

We also obtained general rates
of convergence of the empirical margin distributions to the theoretical
one in the L\'evy distance. Namely, we proved that
the empirical margin distribution converges
to the true margin distribution with probability 1 uniformly over
the class ${\cal G}$ of classifiers if and only if the class ${\cal G}$
is Glivenko-Cantelli. Moreover, if ${\cal G}$ is a Donsker class,
then the rate of convergence in L\'evy distance is $O(n^{-1/4})$.
Faster rates (up to $O(n^{-1/2})$) are possible under some
assumptions on random entropies of the class ${\cal G}.$
We give some examples, showing the optimality of these rates.

We improved previously
known bounds on generalization error of convex combinations
of classifiers.
In particular, our results in Section 3 imply that if the random 
$\eps$-entropy of the
class ${\cal G}$ grows as $\eps^{-\alpha}$ for $\alpha \in (0,2),$
then 
the generalization 
error of any classifier from ${\cal G}$
with zero training error is bounded from above with very high
probability by the quantity
$$
\frac{C}{n^{2/(2+\alpha)}\hat \delta^{2\alpha/(2+\alpha)}},
$$
where $\hat \delta$ is the minimal classification margin of
the training examples and $C$ is a constant. 
The previously known result of Schapire, Freund,  
Bartlett and Lee (1998) gives 
(up to logarithmic factors,  
for ${\cal G}={\rm conv}({\cal H}),$ ${\cal H}$
being a VC-class)  
the bound $O(\frac{1}{n^{1/2}\hat \delta})$ 
which 
corresponds to the worst choice of $\alpha$ ($\alpha=2$). 
We introduce in Section 3 more subtle notions  
of $\gamma$-margin $\delta_n(\gamma ;g)$ and empirical
$\gamma$-margin $\hat \delta_n(\gamma ;g)$ 
(parametrized by $\gamma\in (0,1]$)  
of a classifier $g.$ 
These quantities allow us to obtain similar upper bounds on 
generalization error
of the form 
$$
\frac{C_{\gamma}}{n^{1-\gamma/2}\hat \delta_n^{\gamma}(\gamma ;g)},
$$
in the case when the training error of the classifier $g$ 
is not necessarily equal to $0.$
We call the quantity 
$$
\frac{1}{n^{1-\gamma/2}\hat \delta_n^{\gamma}(\gamma ;g)}
$$
\it the $\gamma$-bound \rm of $g.$   
It follows from the definitions given in Section 3 
that the $\gamma$-bounds decrease when 
$\gamma$ decreases from $1$ to $0.$ 
We prove that for any $\gamma\geq \frac{2\alpha}{2+\alpha}$
with very high probability the $\gamma$-bounds are indeed upper bounds 
on the generalization error (up to a multiplicative constant $C_{\gamma}$).

The proof of the bounds of this type is based
on the powerful concentration inequalities of Talagrand (1996a,b).
For small $\alpha ,$ the bound may become arbitrarily close
to the rate $O(n^{-1}),$ which is known to be the best
possible convergence rate in the zero error case.
In the case of convex combinations of classifiers from
a VC-class ${\cal H},$ one can choose $\alpha =2(V-1)/V,$
where $V$ is the VC-dimension of the class ${\cal H},$
which improves the previously known bounds for convex combinations
of classifiers. We believe that these results can be of importance
in some other learning problems (such as support vector learning,
see Vapnik (1998)).

Koltchinskii, Panchenko and Lozano (2000a,b) studied 
the behavior of the 
$\gamma$-bounds and some other bounds of similar type in 
a number of experiments
with AdaBoost and other methods of combining classifiers. 
We have run AdaBoost for a number 
of rounds with a weak learner that output simple classifiers (e.g. decision stumps)
from a small VC-class. In some of the experiments, we dealt  
with a toy learning problem ("intervals problem") for which 
it was easy to compute the generalization error precisely.
In other cases, we dealt with real data from UCI Irvine 
repository (see Blake and Merz (1998)) and we estimated
the generalization error based on test samples.
In both cases, we computed the $\gamma$-margins and 
the corresponding $\gamma$-bounds based on the training
data and compared the bounds with the generalization error
(or with the test error). 
We give here only a short summary of the results of these experiments
(and some related theoretical results). 
The details are given in Koltchinskii, Panchenko and Lozano (2000a,b). 

$\bullet $ One of the goals of the experiments was to determine the value
of the constant $C_{\gamma}$ involved in the $\gamma$-margin bounds
on generalization error. The results of Section 3 of this paper
show that such a constant exists. Its size, however,
is related to a hard problem of optimizing the constants involved
in Talagrand's concentration inequality for empirical processes 
that was used in the 
proofs. Our experiments showed that the choice $C_{\gamma}=1$ worked  
rather well in the bounds of this type.  
They also showed that the $\gamma$-bounds
did improve the previously known bounds on generalization error 
of AdaBoost. The improvement was significant when the VC-dimension 
of the base class was small and, hence, the parameter $\gamma$ could be 
choosen much smaller than $1.$ Figure 1 shows a typical result
of the experiments.

\begin{figure}
\begin{center}
\epsfig{file=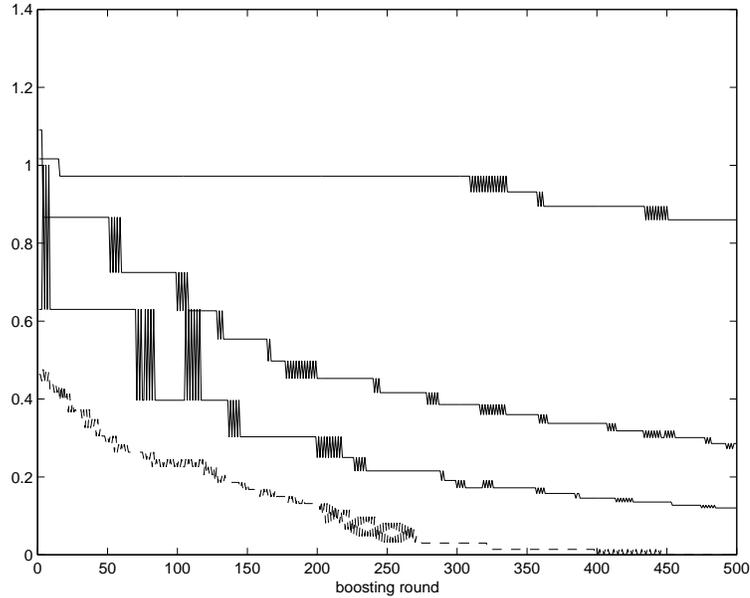,height=8cm,width=10cm}
\caption{Comparison of the generalization error (dashed line) with
the $\gamma$-bounds 
for $\gamma=1,0.8$ and $2/3$ (solid lines, top to bottom)}
\label{figure:intervals1}
\end{center}
\end{figure}

$\bullet$ We also observed that the ratios 
$\frac{\hat \delta_n(\gamma;g)}{\delta_n(\gamma;g)}$
of the empirical $\gamma$-margins to the true $\gamma$-margins of   
classifiers $g$ produced
by AdaBoost had been surprisingly close to $1$ (at least for large sample 
sizes).
The results of Section 3 
imply that, with high probability, 
these ratios are bounded away from $0$ and from $\infty$
uniformly in $g\in {\cal G}$  
for any $\gamma\geq \frac{2\alpha}{2+\alpha}.$ 
Recently, the first 
author proved that the ratios do converge to $1$   
uniformly in $g\in {\cal G}$ a.s. as $n\to \infty$ for  
$\gamma >\frac{2\alpha}{2+\alpha}$ (the example was also given showing that
for $\gamma =\frac{2\alpha}{2+\alpha}$ the ratios do not necessarily converge to $1$ and
for $\gamma <\frac{2\alpha}{2+\alpha}$ they can tend to $\infty$).     
The closeness of the ratios to $1$ explains why the $\gamma$-bounds are 
valid with $C_{\gamma}=1.$

$\bullet$ In the case of the classifiers obtained in consecutive 
rounds of AdaBoost, the $\gamma$-bounds hold even
for the values of $\gamma$ that are substantially smaller
than the threshold $\frac{2\alpha}{2+\alpha}$ given by 
the theory. 
It might be related to the fact that the threshold
is based on the bounds on the entropy of the \it whole \rm
convex hull of the base class ${\cal H}.$ 
On the other hand, AdaBoost 
and other algorithms of this type output classifiers that 
belong to a subset ${\cal G}\subset {\rm conv}({\cal H})$
whose entropy might be much smaller than the entropy of the whole
convex hull.  
Because of this, it is important to develop
adaptive versions of the margin type bounds on generalization
error that take into account the complexity of the classifiers
output by learning algorithms as well as their empirical margins.
A possible approach to this problem was developed in Koltchinskii,
Panchenko and Lozano (2000a).

It should be mentioned that that this paper describes only one
of a number of growing areas of applications of Probability
to computer learning problems. Some other important examples
of such applications are given in Yukich, Stinchcombe and
White (1995), Barron (1991a, b), Barron, Birg\' e and Massart (1999),
Talagrand (1998), Freund (1995, 1999).

\section{Probabilistic bounds for general function classes in terms
of Gaussian and Rademacher complexities}

Let $(S, {\cal A}, P)$ be a probability space and let ${\cal F}$ be a
class of measurable functions from $(S,{\cal A})$ into ${\Bbb R}.$
[Later, in sections 5, 6 we will replace $S$ by $S\times \{-1,1\},$
considering labeled observations; at this point, it is not important].
Let $\{X_k\}$ be a sequence of i.i.d. random variables taking values
in $(S,{\cal A})$ with common distribution $P.$ We assume that this
sequence is defined on a probability space $(\Omega, \Sigma, {\Bbb P}).$
Let $P_n$ be the empirical measure based on the sample $(X_1,\dots ,X_n),$
$$
P_n := n^{-1}\sum_{i=1}^n \delta_{X_i},
$$
where $\delta_x$ denotes the probability distribution concentrated
at the point $x.$
We will denote $Pf := \int_S fdP,$ $P_n f:= \int_S fdP_n,$ etc.

In what follows, $\ell^{\infty}({\cal F})$ denotes the
Banach space of uniformly bounded real valued functions on ${\cal F}$
with the norm
$$
\|Y\|_{\cal F}:= \sup_{f\in {\cal F}}|Y(f)|.
$$

We assume throughout the paper that ${\cal F}$ satisfies standard
measurability assumptions of the theory of empirical processes
(see Dudley (1999), van der Vaart and Wellner (1996))
(for simplicity, one can assume that ${\cal F}$ is countable,
but this, of course, is not necessary).

Our goal in this section is to construct data dependent upper bounds
on the probability $P\{f\leq 0\}$ and on the difference $|P_n\{f\leq 0\}-P\{f\leq 0\}|$
that hold for all $f\in {\cal F}$ with high probability. These inequalities
will be used in the next sections to upper bound the generalization error
of combined classifiers. The bounds will depend on
some measures of "complexity"
of the class ${\cal F}$ which will be introduced next.

Define
$$
G_n({\cal F}):=
{\Bbb E}\|n^{-1}\sum_{i=1}^n g_i\delta_{X_i}\|_{\cal F},
$$
where $\{g_i\}$ is a sequence of i.i.d. standard normal random variables,
independent of $\{X_i\}.$
[Actually, it is common to assume that $\{g_i\}$ is defined on a separate probability
space $(\Omega_g, \Sigma_g, {\Bbb P}_g)$ and that the basic probability space
is now $(\Omega \times \Omega_g, \Sigma \times \Sigma_g, {\Bbb P}\times {\Bbb P}_g)$].
We will call $n\mapsto G_n({\cal F})$ \it the Gaussian complexity function \rm of the
class ${\cal F}$.

Similarly, we define
$$
R_n({\cal F}):=
{\Bbb E}\|n^{-1}\sum_{i=1}^n \eps_i\delta_{X_i}\|_{\cal F},
$$
where $\{\eps_i\}$ is a sequence of i.i.d. Rademacher
(taking values $+1$ and $-1$ with probability $1/2$ each)
random variables,
independent of $\{X_i\}.$
We will call $n\mapsto R_n({\cal F})$ \it
the Rademacher complexity function \rm of the class ${\cal F}$.

One can find in the literature (see, e.g., van der Vaart and Wellner (1996))
various upper bounds on such quantities as $G_n({\cal F})$
and $R_n({\cal F})$ in terms of entropies,
VC-dimensions, etc.


First, we give bounds on $P\{f\leq 0\}$
in terms of a class of so called \it margin cost functions. \rm
These bounds will be used in section 5 in the context of
classification problems to improve recent
results of Mason, Bartlett and Baxter (1999).

Consider a countable family of Lipschitz functions
$\Phi=\{\varphi_k: k\geq 1\},$ where
$\varphi_k:{\Bbb R}\to {\Bbb R}$ are such that
such that $I_{(-\infty,0]}(x)\leq\varphi_k(x)$ for all $k.$
For each $\varphi\in \Phi,\,\,\,$
$L(\varphi)$ will denote its Lipschitz constant.

We assume that for any $x\in S$ the set of real numbers
$\{f(x):f\in {\cal F}\}$ is bounded.

\begin{theorem}
For all $t >0,$
$$
{\Bbb P}\Bigl\{\exists f\in {\cal F}:
P\{f\leq 0\} >
\inf_{k\geq 1}\Bigl[P_n \varphi_k (f)+
4L(\varphi_k)R_n({\cal F})+
\Bigl(\frac{\log k}{n}\Bigr)^{1/2}\Bigr]+
\frac{t}{\sqrt{n}}
\Bigr\}\leq 2\exp \{-2t^2\}
$$
and
$$
{\Bbb P}\Bigl\{\exists f\in {\cal F}:
P\{f\leq 0\} >
\inf_{k\geq 1}\Bigl[P_n \varphi_k (f)+
\sqrt{2\pi}L(\varphi_k)G_n({\cal F})+
\Bigl(\frac{\log k}{n}\Bigr)^{1/2}\Bigr]+
\frac{t+2}{\sqrt{n}}
\Bigr\}\leq 2\exp \{-2t^2\}.
$$
\end{theorem}

{\bf Proof}. Without loss of generality we can and do assume that
each $\varphi \in \Phi$ takes its values in $[0,1]$ (otherwise it
can be redefined as $\varphi \bigwedge 1$). Clearly, in this case
$\varphi (x)=1$ for $x\leq 0.$
For a fixed $\varphi \in \Phi$ and for all $f\in {\cal F}$ we have
\begin{equation}
P\{f\leq 0\} \leq P\varphi (f)\leq P_n \varphi (f)+
\|P_n-P\|_{{\cal G}_{\varphi}},\label{e2.1'}
\end{equation}
where
$$
{\cal G}_{\varphi}:=
\Bigl\{\varphi \circ f -1 : f\in {\cal F}\Bigr\}.
$$
By the exponential inequalities for martingale difference sequences
(see \cite{Dev}, pp 135--136), we have
$$
{\Bbb P}\Bigl\{\|P_n-P\|_{{\cal G}_{\varphi}}\geq
{\Bbb E}\|P_n-P\|_{{\cal G}_{\varphi}} + \frac{t}{\sqrt{n}} \Bigr\}\leq
\exp \{-2t^2\}.
$$
Thus, with probablity at least $1-\exp \{-2t^2\}$ for all $f\in {\cal F}$
\begin{equation}
P\{f\leq 0\} \leq P_n \varphi (f)+
{\Bbb E}\|P_n-P\|_{{\cal G}_{\varphi}} +
\frac{t}{\sqrt{n}}.\label{e2.2'}
\end{equation}
The Symmetrization Inequality gives (\cite{Well})
\begin{equation}
{\Bbb E}\|P_n-P\|_{{\cal G}_{\varphi}}\leq
2{\Bbb E}\|n^{-1}\sum_{i=1}^n \eps_i \delta_{X_i} \|_{{\cal G}_{\varphi}}.
\label{e2.3'}
\end{equation}
Since a function $(\varphi-1)/L(\varphi)$ is a contraction
and $\varphi(0)-1=0,$ the Rademacher comparison inequality
(\cite{LT}, Theorem 4.12, p.112) implies
$$
{\Bbb E}_{\eps}\|n^{-1}\sum_{i=1}^n
\eps_i \delta_{X_i} \|_{{\cal G}_{\varphi}}\leq
2L(\varphi){\Bbb E}_{\eps}\|n^{-1}\sum_{i=1}^n
\eps_i \delta_{X_i} \|_{{\cal F}}.
$$
It now follows from (\ref{e2.2'}), (\ref{e2.3'}) that with probability at least
$1-e^{-2t^2}$ we have for all $f\in {\cal F}$
\begin{equation}
P\{f\leq 0\} \leq P_n \varphi (f)+
4L(\varphi)R_n({\cal F}) +
\frac{t}{\sqrt{n}}.\label{e2.4'}
\end{equation}
We use now (\ref{e2.4'}) with $\varphi=\varphi_k$
and $t$ replaced by $t+\sqrt{\log k}$ to obtain
\begin{eqnarray}
&&
{\Bbb P}\Bigl\{\exists f\in {\cal F}:\
P\{f\leq 0\} > \inf_{k\geq 1}\Bigl[P_n \varphi_k (f)+
4L(\varphi_k)R_n({\cal F}) +
\Bigl(\frac{\log k}{n}\Bigr)^{1/2}\Bigr]+
\frac{t}{\sqrt{n}}\Bigr\}
\nonumber
\\
&&
\leq
\sum_{k\geq 1}\exp\{-2(t+\sqrt{\log {k}})^2\}
\leq
\sum_{k\geq 1}k^{-2}e^{-2t^2}=
\frac{\pi^2}{6}e^{-2t^2}
\leq 2e^{-2t^2}.\label{e2.5'}
\end{eqnarray}

The proof of the second bound is quite similar with the
following changes. The class ${\cal G}_{\varphi}$ is defined
in this case as $\{\varphi \circ f: f\in {\cal F}\}.$
Instead of (\ref{e2.3'}), we have in this case,
by the Symmetrization Inequality and Gaussian Multiplier
Inequality (see \cite{Well}, pp. 108--109,
177--179), that
\begin{equation}
{\Bbb E}\|P_n-P\|_{{\cal G}_{\varphi}}\leq
2{\Bbb E}\|n^{-1}\sum_{i=1}^n \eps_i \delta_{X_i} \|_{{\cal G}_{\varphi}}\leq
\sqrt{2\pi}{\Bbb E}\|n^{-1}\sum_{i=1}^n g_i \delta_{X_i} \|_{{\cal G}_{\varphi}}.
\label{e2.3''}
\end{equation}
Define Gaussian processes
$$
Z_1(f,\sigma):= \sigma n^{-1/2}\sum_{i=1}^n g_i (\varphi \circ f)(X_i)
$$
and
$$
Z_2(f,\sigma):= L(\varphi) n^{-1/2}\sum_{i=1}^n g_i f(X_i)
+\sigma g,
$$
where $\sigma=\pm 1$ and $g$ is standard normal
independent of the sequence $\{g_i\}.$
If we denote by ${\Bbb E}_g$ the expectation on the probability
space $(\Omega_g, \Sigma_g, {\Bbb P}_g)$ on which the sequence
$\{g_i\}$ and $g$ are defined then we have
\begin{equation}
{\Bbb E}_g|Z_1(f,\sigma)-Z_1(h,\sigma')|^2
\leq
{\Bbb E}_g|Z_2(f,\sigma)-Z_2(h,\sigma')|^2,
\label{Add1}
\end{equation}
which is easy to observe if we consider separately the cases
when $\sigma\sigma'$ is equal to $1$ and to $-1.$
Indeed, if $\sigma\sigma'=1$ then (\ref{Add1}) is equivalent to
$$
n^{-1}\sum_{i=1}^n \bigl|\varphi (f(X_i))-\varphi (h(X_i))\bigr|^2
\leq L(\varphi)^2 n^{-1}\sum_{i=1}^n [f(X_i)-h(X_i)]^2
$$
which holds since 
$\varphi$ satisfies the Lipschitz condition with
constant $L(\varphi ).$ 
If $\sigma\sigma'=-1$ then since $0\leq\varphi\leq 1$ we have
\begin{eqnarray*}
&&
{\Bbb E}_g|Z_1(f,\sigma)-Z_1(h,\sigma')|^2 \leq 
2n^{-1}\sum_{i=1}^{n}\varphi^2(f(X_i))+
2n^{-1}\sum_{i=1}^{n}\varphi^2(h(X_i))\leq 
\\
&&
{\Bbb E}(2g)^2\leq
{\Bbb E}_g|Z_2(f,\sigma)-Z_2(h,\sigma')|^2.
\end{eqnarray*}
A version of Slepian's Lemma
(see Ledoux and Talagrand (1991), pp. 76--77) implies that
$$
{\Bbb E}_g\sup \bigl\{Z_1(f,\sigma):
f\in {\cal F},\,\, \sigma=\pm 1\bigr\}
\leq {\Bbb E}_g\sup \bigl\{Z_2(f,\sigma):
f\in {\cal F},\,\, \sigma=\pm 1\bigr\}.
$$
We have
$$
{\Bbb E}_g\|n^{-1/2}\sum_{i=1}^n g_i\delta_{X_i}\|_{{\cal G}_{\varphi}}=
{\Bbb E}_g \sup_{h\in {\bar{\cal G}_{\varphi}}}
\bigl[n^{-1/2}\sum_{i=1}^n g_i h(X_i)\bigr]
={\Bbb E}_g\sup \bigl\{Z_1(f,\sigma):
f\in {\cal F},\sigma=\pm 1\bigr\},
$$
where
${\bar {\cal G}_{\varphi}}:=
\bigl\{\varphi (f), -\varphi(f): f\in {\cal F}
\bigr\},
$
and similarly
$$
L(\varphi){\Bbb E}_g\|n^{-1/2}\sum_{i=1}^n g_i\delta_{X_i}\|_{\cal F}
+{\Bbb E}|g|
\geq
{\Bbb E}_g\sup \bigl\{Z_2(f,\sigma):
f\in {\cal F},\,\,\sigma=\pm 1\bigr\}.
$$
This immediately gives us
\begin{equation}
{\Bbb E}_g\|n^{-1}\sum_{i=1}^n g_i\delta_{X_i}\|_{{\cal G}_{\varphi}}
\leq
L(\varphi)
{\Bbb E}_g\|n^{-1}\sum_{i=1}^n g_i\delta_{X_i}\|_{\cal F}
+n^{-1/2}{\Bbb E}|g|
.
\label{e2.6}
\end{equation}
It follows from (\ref{e2.2'}), (\ref{e2.3''}) and (\ref{e2.6}) that
with probability at least $1-e^{-2t^2}$
\begin{equation}
P\{f\leq 0\} \leq P_n \varphi (f)+
\sqrt{2\pi}L(\varphi)G_n({\cal F}) +
\frac{t+2}{\sqrt{n}}.\label{e2.4''}
\end{equation}
The proof now can be completed the same way as in the case of the
first bound.
\qed

Let us consider a special family of cost functions.
Assume that
$\varphi$ is a fixed \it nonincreasing \rm function
such that $\varphi(x)\geq I_{(-\infty ,0]}(x)$ for $x\in {\Bbb R}$
and $\varphi$ satisfies Lipschitz condition with constant
$L(\varphi).$ Let
$$
\Phi_0 :=\{\varphi(\cdot/\delta) : \delta\in (0,1]\}.
$$
One can easily observe that
$L(\varphi(\cdot/\delta))\leq L(\varphi)\delta^{-1}.$
For this family, Theorem 1 easily implies the following
statement, which, in turn, implies
the result of Schapire, Freund, Bartlett and Lee (1998)
for VC-classes of base classifiers (see Section 5).

\begin{theorem}
For all $t >0,$
\begin{eqnarray*}
{\Bbb P}\Bigl\{\exists f\in {\cal F}:
P\{f\leq 0\}
&>&
\inf_{\delta \in (0,1]}\Bigl[P_n \varphi ({f\over {\delta}})+
{{8L(\varphi)}\over {\delta}}R_n({\cal F})
\\
&+&
\Bigl(\frac{\log\log_2 (2\delta^{-1})}{n}\Bigr)^{1/2}\Bigr]+
\frac{t}{\sqrt{n}}
\Bigr\}
\leq 2\exp \{-2t^2\}
\end{eqnarray*}
and
\begin{eqnarray*}
{\Bbb P}\Bigl\{\exists f\in {\cal F}:
P\{f\leq 0\}
&>&
\inf_{\delta \in (0,1]}\Bigl[P_n \varphi ({f\over {\delta}})+
{{2\sqrt{2\pi}L(\varphi)}\over {\delta}}G_n({\cal F})
\\
&+&
\Bigl(\frac{\log\log_2 (2\delta^{-1})}{n}\Bigr)^{1/2}\Bigr]+
\frac{t+2}{\sqrt{n}}
\Bigr\}
\leq 2\exp \{-2t^2\}.
\end{eqnarray*}
\end{theorem}

{\bf Proof.}
One has to apply the bounds of Theorem 1
for the sequence $\varphi_k(\cdot):=\varphi(\cdot/\delta_k),$
where $\delta_k=2^{-k},$
and then notice that for $\delta \in (\delta_k, \delta_{k-1}],$
we have
$$
\frac{1}{\delta_k}\leq \frac{2}{\delta},\
P_n\varphi({f\over {\delta_k}})\leq P_n\varphi ({f\over {\delta}})
$$
and
$$
\sqrt{\log k}=\sqrt{\log \log_2{\frac{1}{\delta_k}}}\leq
\sqrt{\log \log_2{\frac{2}{\delta}}}.
$$
\qed

{\bf Remark}. The constant $8$ in front of the Rademacher complexity
and the constant $2\sqrt{2\pi}$ in front of the Gaussian complexity
can be replaced by $4c$ and $\sqrt{2\pi}c,$ respectively, for any
$c>1$ (with minor changes in the logarithmic term). Also, one can
choose $c=c(\delta),$ where $c(\delta)=1+o(1)$ as $\delta \to 0.$
 
In the next statements we use the Rademacher complexities, but
Gaussian complexities can be used similarly.

Assuming now that
$\varphi$ is a function from ${\Bbb R}$ into ${\Bbb R}$ such that
$\varphi (x)\leq I_{(-\infty ,0]}(x)$ for all $x\in {\Bbb R}$ and
$\varphi $ still satisfies the Lipschitz condition with constant
$L(\varphi),$
one can prove the following statement.

\begin{theorem}
For all $t >0,$
\begin{eqnarray*}
{\Bbb P}\Bigl\{\exists f\in {\cal F}:
P\{f\leq 0\}
&<&
\sup_{\delta \in (0,1]}
\Bigl(P_n \varphi ({f\over {\delta}})-
{{8L(\varphi)}\over {\delta}}R_n({\cal F})
\\
&-&
\Bigl(\frac{\log\log_2 (2\delta^{-1})}{n}\Bigr)^{1/2}\Bigr)-
\frac{t}{\sqrt{n}}
\Bigr\}
\leq 2\exp \{-2t^2\}.
\end{eqnarray*}
\end{theorem}

Denote
$$
\Delta_n ({\cal F};\delta):=
{8\over {\delta}}R_n({\cal F})+
\Bigl(\frac{\log\log_2 (2\delta^{-1})}{n}\Bigr)^{1/2}.
$$

The bounds of theorems 2 and 3 easily imply that for all $t>0$
$$
{\Bbb P}\Bigl\{\exists f\in {\cal F}:
P\{f\leq 0\} > P_n\{f\leq 0\} +
\inf_{\delta \in (0,1]}\Bigl[P_n\{0<f\leq \delta\} +
\Delta_n({\cal F};\delta)\Bigr]+\frac{t}{\sqrt{n}}
\Bigr\}
\leq 2\exp \{-2t^2\}
$$
and
$$
{\Bbb P}\Bigl\{\exists f\in {\cal F}:
P\{f\leq 0\} < P_n\{f\leq 0\} -
\inf_{\delta \in (0,1]}\Bigl[P_n\{-\delta<f\leq 0\} +
\Delta_n({\cal F};\delta)
\Bigr]-\frac{t}{\sqrt{n}}
\Bigr\}
\leq 2\exp \{-2t^2\}.
$$
To prove this it's enough to take $\varphi$ equal to $1$ for $x\leq 0,$
$0$ for $x\geq 1$ and linear in between in the case of the first bound;
in the case of the second bound, the choice of $\varphi$ is
$1$ for $x\leq -1,$ $0$ for $x\geq 0$ and linear in between.
Similarly, it can be shown that
$$
{\Bbb P}\Bigl\{\exists f\in {\cal F}:
P_n\{f\leq 0\} > P\{f\leq 0\} +
\inf_{\delta \in (0,1]}\Bigl[P\{0<f\leq \delta\} +
\Delta_n({\cal F};\delta)
\Bigr]+\frac{t}{\sqrt{n}}
\Bigr\}
\leq 2\exp \{-2t^2\}
$$
and
$$
{\Bbb P}\Bigl\{\exists f\in {\cal F}:
P_n\{f\leq 0\} < P\{f\leq 0\} -
\inf_{\delta \in (0,1]}\Bigl[P\{-\delta<f\leq 0\} +
\Delta_n({\cal F};\delta)
\Bigr]-\frac{t}{\sqrt{n}}
\Bigr\}
\leq 2\exp \{-2t^2\}.
$$

Combining the last bounds, we get the following result:

\begin{theorem}
For all $t >0,$
$$
{\Bbb P}\Bigl\{\exists f\in {\cal F}:
|P_n\{f\leq 0\}-P\{f\leq 0\}|>
\inf_{\delta \in (0,1]}\Bigl[P_n\{|f|\leq \delta\}+
\Delta_n({\cal F};\delta)
\Bigr]+\frac{t}{\sqrt{n}}
\Bigr\}
\leq 4\exp \{-2t^2\}
$$
and
$$
{\Bbb P}\Bigl\{\exists f\in {\cal F}:
|P_n\{f\leq 0\}-P\{f\leq 0\}|>
\inf_{\delta \in (0,1]}\Bigl[P\{|f|\leq \delta\}+
\Delta_n({\cal F};\delta)
\Bigr]+\frac{t}{\sqrt{n}}
\Bigr\}
\leq 4\exp \{-2t^2\}.
$$
\end{theorem}

Denote
$$
H_f(\delta):=\delta P\{|f|\leq \delta\},\
H_{n,f}(\delta):=\delta P_n\{|f|\leq \delta\}.
$$
Plugging in the second bound of Theorem 4
$\delta:=H_f^{-1}(R_n({\cal F}))\bigwedge 1$ 
(we use the notation $a \bigwedge b:=\min(a,b)$)
easily gives
us the following upper bound that holds
for any $t>0$ with probability at least $1-4e^{-2t^2}:$
$$
\forall f\in {\cal F}\
 |P_n\{f\leq 0\}-P\{f\leq 0\}|\leq
\frac{9R_n({\cal F})}{\delta}+
\Bigl(\frac{\log\log_2 (2\delta^{-1})}{n}\Bigr)^{1/2}
+\frac{t}{\sqrt{n}}.
$$
Similarly, the first bound of Theorem 4 gives that
for any $t>0$ with probability
at least $1-4e^{-2t^2}:$
$$
\forall f\in {\cal F}\ |P_n\{f\leq 0\}-P\{f\leq 0\}|\leq
\frac{9R_n({\cal F})}{\delta}+
\Bigl(\frac{\log\log_2 (2\delta^{-1})}{n}\Bigr)^{1/2}
+\frac{t}{\sqrt{n}}
$$
with $\delta:=H_{n,f}^{-1}(R_n({\cal F}))\bigwedge 1.$


The next example shows that, in general, the term
${1\over{\delta}}R_n({\cal F})$
of the bound of Theorem 2 (and other similar results, in particular,
Theorem 4) can not be improved.

Let us consider a sequence
$\{X_n\}$
of independent identically distributed random variables in
$l_{\infty}$
defined by
$$
X_n=\left\{
\eps_k^n (2\log(k+1))^{-\frac{1}{2}}
\right\}_{k\geq 1},\ n\geq 1,
$$
where
$\eps_k^n$
are i.i.d. Rademacher random variables
($P(\eps_k^n=\pm 1)=1/2$).
We consider a class of functions that consists
of canonical projections on each coordinate
$$
{\cal F}=\{f_k: f_k(x)=x_k\}.
$$

Let
$\phi(x)$
be an increasing function such that
$\phi(0)=0.$
Then the following proposition holds.

\begin{proposition}
$$
{\Bbb P}\Biggl\{
\exists f\in {\cal F}: P\{f\leq 0\}\geq
\inf_{\delta\in (0,1]}[
P_n\{f\leq \delta\}+
\frac{8}{\phi(\delta)}
R_n({\cal F})]+\frac{t}{\sqrt{n}}
\Biggr\}\to 1
$$
when
$n\to \infty$
uniformly for all
$t\leq 2^{-1}n^{1/2}\phi((4n)^{-1/2})-c,$
where $c>0$ is some fixed constant.
\end{proposition}

\vspace{2mm}

{\bf Proof.}
It's well known that
${\cal F}$
is a bounded CLT class for the distribution
$P$
of the sequence
$\{X_n\}$
(see Ledoux and Talagrand (1991), pp. 276--277).
Notice that
$P(f_k\leq 0)=1/2$
for all $k$
and
${\Bbb E}\|n^{-1}\sum \eps_i\delta_{X_i}\|_{\cal F}\leq cn^{-1/2}$
for some constant
$c>0.$
Let us denote by
$t'=t+2\sqrt{2\pi}c.$
The infimum inside the probability is less then or equal to
the value of the expression at any fixed point. Therefore,
for each $k$ we will choose
$\delta$
to be equal to a
$\delta_k>(2\log(k+1))^{-1/2}.$
It's easy to see that for this value of
$\delta,$
$$
P_n\{f_k\leq \delta_k\}
=\frac{1}{n}\sum_{i=1}^n I(\eps_k^i=-1).
$$
Combining these estimates we get that
the probability defined in the statement of the proposition
is greater than or equal to
\begin{eqnarray*}
{\Bbb P}\Biggl\{
\exists k:
\frac{1}{2}\geq \frac{1}{n}\sum_{i\leq n}I(\eps_k^i =-1)+
\frac{t'}{\phi(\delta_k)\sqrt{n}}
\Biggr\}
= 1-\prod_{k}
{\Bbb P}\Biggl\{
\frac{1}{2}< \frac{1}{n}\sum_{i\leq n}I(\eps_k^i =-1)+
\frac{t'}{\phi(\delta_k)\sqrt{n}}
\Biggr\}
\end{eqnarray*}
In the product above factors are possibly not equal to
$1$
only for $k$ in the set of indices
$$
{\cal K}=\Bigl\{k:\gamma_k=
\frac{t'}{\phi(\delta_k)\sqrt{n}}\leq \frac{1}{2}
\Bigr\}.
$$
Clearly,
$$
{\Bbb P}\Biggl\{
1/2< n^{-1}\sum_{i\leq n}I(\eps_1^i = -1)+\delta
\Biggr\}\leq
1-{n\choose k_0}2^{-n},
$$
where
$k_0=[n/2-\delta n]-1.$
For simplicity of calculations we will set
$k_0=n/2-\delta n.$
Utilizing the following estimates in Stirling's formula
for the factorial (see Feller (1950))
\begin{equation}
(2\pi)^{\frac{1}{2}}n^{n+\frac{1}{2}}e^{-n+1/(12n+1)}<
n!<
(2\pi)^{\frac{1}{2}}n^{n+\frac{1}{2}}e^{-n+1/12n}\label{Stirling}
\end{equation}
it is straightforward to check that
for some constant $c>0$
\begin{equation}
{n\choose k_0}2^{-n}\geq cn^{-\frac{1}{2}}\left(
(1-2\delta)^{1-2\delta}(1+2\delta)^{1+2\delta}
\right)^{-\frac{n}{2}}\geq
cn^{-\frac{1}{2}}\exp(-4n\delta^2).\label{Lower}
\end{equation}
The last inequality is due to the fact that
$$
\exp(x^2)\leq
(1-x)^{1-x}(1+x)^{1+x}\leq \exp(2x^2)\label{expbounds}
$$
for $x<2^{-1/2}.$
It follows from (\ref{Lower}) that
$$
{\Bbb P}\Biggl\{
\frac{1}{2}< \frac{1}{n}\sum_{i\leq n}I(\eps_k^i =-1)+
\gamma_k
\Biggr\}\leq 1-cn^{-1/2}\exp(-4n\gamma_k^2).
$$
Since $\gamma_k\leq 1/2$ for $k\in {\cal K},$
we can continue and come to the following lower bound
\begin{eqnarray*}
&&
1-\prod_{k\in {\cal K}}(1-cn^{-1/2}\exp(-4n\gamma_k^2))\geq
1-\exp(-\sum_{k\in {\cal K}}cn^{-1/2}\exp(-4n\gamma_k^2))
\\
&&
\geq
1-\exp(-{\rm card}({\cal K})cn^{-1/2}e^{-n})\to 1,
\end{eqnarray*}
uniformly in $t',$ if we check that
${\rm card}({\cal K})cn^{-1/2}e^{-n}\to \infty.$
Indeed, if
$$
t'\leq 2^{-1}n^{1/2}\phi((4n)^{-1/2})
$$
then for $n$ large enough
$$
t'\leq 2^{-1}n^{1/2}\phi((4n)^{-1/2})\leq
2^{-1}n^{1/2}\phi((2\log([cne^{n}]+1))^{-1/2}).
$$
It means that
$[cne^n]\in {\cal K},$ and, therefore,
$${\rm card}({\cal K})cn^{-1/2}e^{-n}\geq n^{1/2}-{1\over {cn^{1/2}e^n}}\to \infty .$$
Proposition is proven.

\qed
\vspace{2mm}

{\bf Remarks.} If
$\phi(x)=x^{1-\alpha}$ for some positive
$\alpha$
then the convergence in the proposition
holds for
$t\leq cn^{\alpha/2}.$
Also, if ${{\phi (\delta)}\over {\delta}}\to \infty$ as $\delta \to 0,$
then the convergence in the proposition holds uniformly in $t\in [0,T]$
for any $T>0.$ It means that the bound of Theorem 2 does not hold with
$\frac{1}{\delta}R_n({\cal F})$ replaced by
$\frac{1}{\phi(\delta)}R_n({\cal F})$.
Similarly, one can show that
$$
{\Bbb P}\Biggl\{
\exists f\in {\cal F}: |P_n\{f\leq 0\}-P\{f\leq 0\}|\geq
\inf_{\delta\in (0,1]}[
P_n\{|f|\leq \delta\}+
\frac{8}{\phi(\delta)}
R_n({\cal F})]+\frac{t}{\sqrt{n}}
\Biggr\}\to 1
$$
when
$n\to \infty$
uniformly for all
$t\leq 2^{-1}n^{1/2}\phi((4n)^{-1/2})-c.$

\mdsk


\mdsk

\section{Conditions on random entropies and $\gamma$-margins}

Given a metric space $(T,d),$ we denote $H_d(T;\eps)$
the $\eps$-entropy of $T$ with respect to $d,$ i.e.
$$
H_d(T;\eps):=\log N_d(T;\eps),
$$
where $N_d(T;\eps)$ is the minimal number of balls of
radius $\eps$ covering $T.$
Let $d_{P_n,2}$ denote the metric of the space $L_2(S;dP_n):$
$$
d_{P_n,2}(f,g):=\bigl(P_n |f-g|^2\bigr)^{1/2}.
$$

The next theorems improve the bounds of previous section
under some assumptions on the growth of random entropies
$H_{d_{P_n,2}}({\cal F};\cdot).$ We will use these results in
section 5 to obtain an improvement of the bound of Schapire,
Freund, Bartlett and Lee (1998) on generalization error
of boosting. The method of proof is similar to the one
developed in Koltchinskii and Panchenko (1999) and is
based on powerful concentration inequalities of Talagrand (1996)
(see also Massart (2000)).

Define for $\gamma\in (0,1]$
$$
\delta_n(\gamma ;f):=
\sup\Bigl\{\delta\in (0,1):
\delta^{\gamma}
P\{f\leq \delta\}\leq n^{-1+\frac{\gamma}{2}}
\Bigr\}
$$
and
$$
\hat \delta_n(\gamma ;f):=
\sup\Bigl\{\delta\in (0,1):
\delta^{\gamma }
P_n\{f\leq \delta\}\leq n^{-1+\frac{\gamma}{2}}
\Bigr\}.
$$
We call $\delta_n(\gamma ;f)$ and $\hat \delta_n(\gamma ;f),$
respectively, the \it $\gamma$-margin \rm and the \it empirical
$\gamma$-margin \rm of $f.$

The main result of this section is Theorem 5 that gives the
condition on the random entropy $H_{d_{P_n,2}}({\cal F};\cdot)$
under which the true $\gamma$-margin of any $f\in {\cal F}$
is with probability very close to $1$ within a multiplicative
constant from its empirical $\gamma$-margin. This implies
that with high probability for all $f\in {\cal F}$
$$
P\{f\leq 0\}\leq \frac{\rm const}
{n^{1-\gamma/2}\hat \delta_n(\gamma ;f)^{\gamma}}.
$$
The bounds of previous section correspond to the case of $\gamma =1.$
It is easy to see from the definitions of $\gamma-$margins that
the quantity
$(n^{1-\gamma/2}\hat \delta_n (\gamma ;f)^{\gamma})^{-1}$
(called in the introduction the $\gamma$-bound) 
increases in $\gamma \in (0,1].$ This shows that the bound
in the case of $\gamma <1$ is tighter than the bounds of Section 2.

\mdsk

\begin{theorem}
Suppose that for some $\alpha \in (0,2)$ and for some constant $D>0$
\begin{equation}
H_{d_{P_n,2}}\bigl({\cal F};u\bigr)\leq Du^{-\alpha},\ u>0\ {\rm a.s.}
\label{entropy}
\end{equation}
Then for any $\gamma \geq \frac{2\alpha}{2+\alpha},$ for
some constants $A, B>0$ and for all large enough $n$
\begin{eqnarray*}
{\Bbb P}
\Bigl\{\forall f\in {\cal F}:
A^{-1}\hat \delta_n(\gamma ;f)
\leq
\delta_n({\gamma ;f})\leq
A\hat \delta_n(\gamma ;f)\Bigr\}
\geq
1-B\log_2\log_2{n}
\exp\Bigl\{-n^{\frac{\gamma}{2}}/2 \Bigr\}.
\end{eqnarray*}
\end{theorem}

The proof is based on the following result.

\begin{theorem}
Suppose that for some $\alpha \in (0,2)$ and for some
constant $D>0$
condition (\ref{entropy}) holds.
Then for some constants $A, B>0,$ for all
$\delta \geq 0$
and
\begin{equation}
\varepsilon \geq
\Bigl(\frac{1}{n\delta^{\alpha}}\Bigr)^{\frac{2}{2+\alpha}}
\vee \frac{2\log n}{n},
\label{condeps}
\end{equation}
and for all large enough $n,$
the following bounds hold:
$$
{\Bbb P}\Bigl\{\exists f\in {\cal F}\
P_n\{f\leq \delta\}\leq
\varepsilon\
{\rm and}\
P\{f\leq  \frac{\delta}{2}\}\geq A
\varepsilon
\Bigr\}
\leq B \log_2\log_2\varepsilon^{-1}
\exp\{-\frac{n\varepsilon}{2}\}.
$$
and
$$
{\Bbb P}\Bigl\{\exists f\in {\cal F}\
P\{f\leq \delta\}\leq
\varepsilon\
{\rm and}\
P_n\{f\leq  \frac{\delta}{2}\}\geq A\varepsilon
\Bigr\}
\leq B \log_2\log_2\varepsilon^{-1}
\exp\{-\frac{n\varepsilon}{2}\}.
$$

\end{theorem}

\mdsk

{\bf Proof}.
Define recursively
$$
r_0:=1,\
r_{k+1}=
C\sqrt{r_k \varepsilon}\bigwedge 1
$$
with some sufficiently large constant $C>1$ (the choice of $C$
will be explained later).
By a simple induction argument we have either
$C\sqrt{\eps}\geq 1$ and $r_k\equiv 1,$ or
$C\sqrt{\eps}< 1$ and in this case
$$
r_k=C^{1+2^{-1}+\dots +2^{-(k-1)}}\eps^{2^{-1}+\dots +2^{-k}}=
C^{2(1-2^{-k})}\eps^{1-2^{-k}}= (C\sqrt{\eps})^{2(1-2^{-k})}.
$$
Without loss of generality we can assume that 
$C\sqrt{\eps}< 1.$
Let
$$
\gamma_k := \sqrt{\frac{\eps}{r_k}}=C^{2^{-k}-1}\eps^{2^{-k-1}}.
$$
For a fixed $\delta > 0,$ define
$$
\delta_0 = \delta,\,
\delta_k := \delta (1-\gamma_0-\dots -\gamma_{k-1}),\,
\delta_{k,\frac{1}{2}}=\frac{1}{2}(\delta_k + \delta_{k+1}),\,
\ k\geq 1.
$$

{\bf Warning.} In what follows in the proof
``$c$'' denotes a constant; its values can be different in different
places.

Define ${\cal F}_0:={\cal F},$ and further recursively
$$
{\cal F}_{k+1}:=\Bigl\{f\in {\cal F}_k:
P\{f\leq \delta_{k,\frac{1}{2}}\}\leq {r_{k+1}/2}\Bigr\}.
$$
For $k\geq 0,$
let $\varphi_k $ be a continuous function from ${\Bbb R}$ into $[0,1]$
such that $\varphi_k (u)=1$ for $u\leq \delta_{k,\frac{1}{2}},$
$\varphi_k (u)=0$ for $u\geq \delta_k,$  and linear
for $\delta_{k,\frac{1}{2}}\leq u \leq \delta_k.$
For $k\geq 1$
let $\varphi_k^{\prime} $ be a continuous function from ${\Bbb R}$ into $[0,1]$
such that $\varphi_k^{\prime} (u)=1$ for $u\leq \delta_k,$
$\varphi_k^{\prime} (u)=0$ for $u\geq \delta_{k-1,\frac{1}{2}},$
and linear for $\delta_k\leq u \leq
\delta_{k-1,\frac{1}{2}}.$
We have
\begin{eqnarray*}
&&
\sum_{i=0}^{k}\gamma_{i}=
C^{-1}\bigl[C\sqrt{\eps}+(C\sqrt{\eps})^{2^{-1}}+\dots
+(C\sqrt{\eps})^{2^{-k}}\bigr]
\\
&&
\leq
C^{-1}(C\sqrt{\eps})^{2^{-k}}
(1-(C\sqrt{\eps})^{2^{-k}})^{-1}\leq 1/2,
\end{eqnarray*}
for $\eps\leq C^{-4},$ $C>2(2^{1/4}-1)^{-1}$ and
$k\leq \log_2\log_2 \eps^{-1}.$
Hence, for small enough $\eps$
(note that our choice of
$\eps\leq C^{-4}$ implies $C\sqrt{\eps}<1$), we have
$$
\gamma_0+\dots +\gamma_{k} \leq \frac{1}{2},\ k\geq 1.
$$
Therefore, for all $k\geq 1,$ we get $\delta_k \in (\delta/2, \delta).$
Note also that below our choice of $k$
will be such that the restriction
$k\leq \log_2\log_2 \eps^{-1}$
for any fixed $\eps>0$ will always be fulfilled.

Define
$$
{\cal G}_k :=
\bigl\{\varphi_k\circ f:f\in {\cal F}_k\bigr\},
\,\,\,k\geq 0
$$
and
$$
{\cal G}_k^{\prime} :=
\bigl\{\varphi_k^{\prime}\circ f:f\in {\cal F}_{k}\bigr\},
\,\,\,k\geq 1 .
$$
Clearly, by these definitions, for $k\geq 1$
$$
\sup_{g\in {\cal G}_k} Pg^2 \leq
\sup_{f\in {\cal F}_k}P\{f\leq \delta_k\}\leq
\sup_{f\in {\cal F}_k}P\{f\leq \delta_{k-1,\frac{1}{2}}\}
\leq r_k/2 \leq r_k
$$
and
$$
\sup_{g\in {\cal G}_k^{\prime}} Pg^2 \leq
\sup_{f\in {\cal F}_{k}}P\{f\leq \delta_{k-1,\frac{1}{2}}\}
\leq r_{k}/2\leq r_k .
$$
Since $r_0=1,$ for $k=0$ the first
inequality becomes trivial.
If now we introduce the following events
$$
E^{(k)} := \Bigl\{\|P_n-P\|_{{\cal G}_{k-1}}\leq
K_1 {\Bbb E}\|P_n-P\|_{{\cal G}_{k-1}}+ K_2\sqrt{r_{k-1}\eps}+K_3\eps\Bigr\}
\bigcap
$$
$$
\bigcap \Bigl\{\|P_n-P\|_{{\cal G}_k^{\prime}}\leq
K_1 {\Bbb E}\|P_n-P\|_{{\cal G}_k^{\prime}}+ K_2\sqrt{r_{k}\eps}+
K_3\eps\Bigr\},\,\,k\geq 1,
$$
then it follows from the concentration inequalities
of Talagrand (1996a,b)
(see also \cite{Mas1}) that with some numerical constants
$K_1, K_2, K_3>0$
$$
{\Bbb P}((E^{(k)})^c)\leq 2e^{-\frac{n\eps}{2}}.
$$
Denote $E_0=\Omega,$
$$
E_N:=\bigcap_{k=1}^N E^{(k)},
\,\,\,N\geq 1.
$$
Then
$$
{\Bbb P}(E_N^c)\leq 2N e^{-\frac{n\eps}{2}}.
$$
In what follows we can and do assume without loss of generality
that $\eps <C^{-4}$ and therefore, $r_{k+1}<r_k$
and $\delta_k\in (\delta/2,\delta],\
k\leq \log_2\log_2 \eps^{-1}.$
(If $\eps\geq C^{-4},$ then the bounds of the theorem
obviously hold with any constant $A>C^4.$)
The following lemma holds.

\begin{lemma}
Let $N$ be such that
\begin{equation}
N\leq \log_2\log_2 \eps^{-1}\,\,\,
\mbox{ and }\,\,\,
r_N\geq \eps
\label{mystars}.
\end{equation}
Let 
${\cal J} =
\Bigl\{\inf_{f\in {\cal F}}P_n\{f\leq \delta\}\leq \eps\Bigr\}.$
Then the following properties hold on the event $E_N\bigcap {\cal J}:$
$$
(i)\ \forall f\in {\cal F}\ P_n\{f\leq \delta\}\leq \eps \Longrightarrow
f\in {\cal F}_N
$$
and
$$
(ii)\ \sup_{f\in {\cal F}_k}P_n\{f\leq \delta_k\}\leq r_k,\
0\leq k\leq N.
$$

\end{lemma}

\mdsk

{\bf Proof}. We will use the induction with respect to $N.$
For $N=0,$ the statement is obvious. Suppose it holds for
some $N\geq 0,$ such that $N+1$ still satisfies condition
(\ref{mystars}) of the lemma.
Then on the event $E_{N}\bigcap {\cal J}$ we have
$$
\sup_{f\in {\cal F}_k}P_n\{f\leq \delta_k\}\leq r_k,\
0\leq k\leq N
$$
and
$$
\forall f\in {\cal F}\ P_n\{f\leq \delta\}\leq \eps \Longrightarrow
f\in {\cal F}_N.
$$
Suppose now that $f\in {\cal F}$ is such that
$P_n\{f\leq \delta\}\leq \eps. $ By the induction assumptions,
on the event $E_N,$ we have $f\in {\cal F}_N.$
Because of this, we obtain on the event $E_{N+1}$
\begin{eqnarray}
&&
P\{f\leq \delta_{N,\frac{1}{2}}\}
\leq P_n\{f\leq \delta_N\}+\|P_n-P\|_{{\cal G}_N}
\nonumber
\\
&&
\leq \eps
+K_1 {\Bbb E}\|P_n-P\|_{{\cal G}_{N}}+ K_2\sqrt{r_{N}\eps}+
K_3\eps  .\label{e3.1}
\end{eqnarray}
For a class ${\cal G},$ define
$$
\hat R_n ({\cal G}):=\|n^{-1}\sum_{i=1}^n \eps_i \delta_{X_i}\|_{\cal G},
$$
where $\{\eps_i\}$ is a sequence
of i.i.d. Rademacher random variables.
By the symmetrization inequality,
\begin{equation}
{\Bbb E}\|P_n-P\|_{{\cal G}_{N}}\leq
2{\Bbb E}I_{E_N} {\Bbb E}_{\eps} \hat R_n({\cal G}_{N})+
2{\Bbb E}I_{E_N^c}{\Bbb E}_{\eps}\hat R_n({\cal G}_{N}).
\label{e3.2}
\end{equation}
Next, by the well known entropy inequalities for subgaussian processes
(see van der Vaart and Wellner (1996), Corollary 2.2.8), we have
\begin{equation}
{\Bbb E}_{\eps} \hat R_n({\cal G}_{N})\leq
\inf_{g\in {\cal G}_{N}}
{\Bbb E}_{\eps}\bigl|n^{-1}\sum_{j=1}^n \eps_j g(X_j)\bigr|
+\frac{c}{\sqrt{n}}
\int_0^{(2\sup_{g\in {\cal G}_{N}}P_n g^2)^{1/2}}
H_{d_{P_n,2}}^{1/2}({\cal G}_{N};u)du.
\label{e3.3}
\end{equation}
By the induction assumption, on the event $E_N\bigcap {\cal J}$
\begin{eqnarray*}
&&
\inf_{g\in {\cal G}_{N}}
{\Bbb E}_{\eps}\bigl|n^{-1}\sum_{j=1}^n \eps_j g(X_j)\bigr|\leq
\inf_{g\in {\cal G}_{N}}
{\Bbb E}_{\eps}^{1/2}\bigl|n^{-1}\sum_{j=1}^n \eps_j g(X_j)\bigr|^2\leq
\frac{1}{\sqrt{n}}\inf_{g\in {\cal G}_{N}}\sqrt{P_n g^2}
\\
&&
\leq \frac{1}{\sqrt{n}}
\inf_{f\in {\cal F}_{N}}\sqrt{P_n \{f\leq \delta_N\}}\leq
\frac{1}{\sqrt{n}}
\inf_{f\in {\cal F}_{N}}\sqrt{P_n \{f\leq \delta\}}\leq
\sqrt{\frac{\eps}{n}}\leq \eps
\end{eqnarray*}
We also have on the event $E_N\bigcap {\cal J}$
$$
\sup_{g\in {\cal G}_{N}}P_n g^2\leq \sup_{f\in {\cal F}_N}
P_n\{f\leq \delta_N\}\leq r_N.
$$
The Lipschitz norm of $\varphi_{k-1}$ and $\varphi_k^{\prime}$ is
bounded by
\begin{eqnarray*}
L=2(\delta_{k-1}-\delta_k)^{-1}=
2\delta^{-1}\gamma_{k-1}^{-1}=
\frac{2}{\delta}\sqrt{\frac{r_{k-1}}{\eps}}
\end{eqnarray*}
which implies the following bound on the distance
$$
d_{P_n,2}^2\Bigl(\varphi_N \circ f;\varphi_N\circ g\Bigr)
=
n^{-1}\sum_{j=1}^n \Bigl|\varphi_N (f(X_j))-
\varphi_N (g(X_j))\Bigr|^2\leq
\Bigl(\frac{2}{\delta}\sqrt{\frac{r_{N}}{\eps}}\Bigr)^2 
d_{P_n,2}^2 (f,g).
$$
Therefore, on the event $E_N\bigcap {\cal J}$
\begin{eqnarray}
&&
\frac{1}{\sqrt{n}}
\int_0^{(2\sup_{g\in {\cal G}_{N}}P_n g^2)^{1/2}}
H_{d_{P_n,2}}^{1/2}({\cal G}_{N};u)du
\leq
\frac{1}{\sqrt{n}}
\int_0^{(2r_N)^{1/2}}
H_{d_{P_n,2}}^{1/2}({\cal F};
\frac{\delta \sqrt{\eps} u}{2\sqrt{r_N}})du
\nonumber
\\
&&
\leq
c(\frac{r_N}{\eps})^{\alpha /4}
\frac{r_N^{1/2-\alpha/4}}{\sqrt{n}\delta^{\alpha/2}}
\leq c\frac{r_N^{1/2}}{\eps^{\alpha/4}}\eps^{\frac{2+\alpha}{4}}
= c\sqrt{r_N \eps},
\label{e3.6}
\end{eqnarray}
where we used the fact that condition (\ref{condeps}) of the theorem
implies
$$
\frac{1}{n^{1/2}\delta^{\alpha/2}}\leq \eps^{\frac{2+\alpha}{4}}.
$$
It follows from (\ref{e3.3}), (\ref{e3.6}) that on the event
$E_{N+1} \bigcap {\cal J}$
\begin{equation}
{\Bbb E}_{\eps} \hat R_n({\cal G}_{N})\leq
c\sqrt{r_N \eps}.\label{e3.7}
\end{equation}
Since we also have
$$
{\Bbb E}_{\eps} \hat R_n({\cal G}_{N+1})\leq 1,
$$
(\ref{e3.2}) and (\ref{e3.7}) yield
$$
{\Bbb E}\|P_n-P\|_{{\cal G}_{N}}\leq
c\sqrt{r_N \eps}+
2{\Bbb P}(E_N^c)
\leq
c\sqrt{r_N \eps}+
4Ne^{-n\eps/2}.
$$
Since $4Ne^{-n\eps/2}\leq \eps$ (it holds due to the
conditions (\ref{condeps}) and (\ref{mystars}), for all large enough $n$)
we conclude that
with some constant $c>0$
$$
{\Bbb E}\|P_n-P\|_{{\cal G}_{N}}\leq
c\sqrt{r_N \eps}.
$$
Now we use (\ref{e3.1}) and see that on the event
$E_{N+1}\bigcap {\cal J}$
\begin{equation}
P\{f\leq \delta_{N,\frac{1}{2}}\}\leq
c\bigl(\eps+
\sqrt{r_N \eps}\bigr).
\end{equation}
Therefore, it follows that with a proper choice of
constant $C>0$ in the recurrence relationship defining
the sequence $\{r_k\},$ we have on the event
$E_{N+1}\bigcap {\cal J}$
$$
P\{f\leq \delta_{N,\frac{1}{2}}\}\leq
\frac{1}{2}C\sqrt{r_N \eps}=r_{N+1}/2.
$$
This means that $f\in {\cal F}_{N+1}$
and the induction step for (i) is proved.
This will now imply (ii).
We have on the event $E_{N+1}$
\begin{eqnarray}
&&
\sup_{f\in {\cal F}_{N+1}} P_n\{f\leq \delta_{N+1}\}\leq
\sup_{f\in {\cal F}_{N+1}} P\{f\leq \delta_{N,\frac{1}{2}}\}+
\|P_n-P\|_{{\cal G}_{N+1}^{\prime}}
\nonumber
\\
&&
\leq r_{N+1}/2 +
K_1 {\Bbb E}\|P_n-P\|_{{\cal G}_{N+1}^{\prime}}+ K_2\sqrt{r_{N+1}\eps}+
K_3\eps .\label{ee3.1}
\end{eqnarray}
By the symmetrization inequality,
\begin{equation}
{\Bbb E}\|P_n-P\|_{{\cal G}_{N+1}^{\prime}}\leq
2{\Bbb E}I_{E_N} {\Bbb E}_{\eps} \hat R_n({\cal G}_{N+1}^{\prime})+
2{\Bbb E}I_{E_N^c}{\Bbb E}_{\eps}\hat R_n({\cal G}_{N+1}^{\prime}).
\label{ee3.2}
\end{equation}
As above, we have
\begin{equation}
{\Bbb E}_{\eps} R_n({\cal G}_{N+1}^{\prime})\leq
\inf_{g\in {\cal G}_{N+1}^{\prime}}
{\Bbb E}_{\eps}\bigl|n^{-1}\sum_{j=1}^n \eps_j g(X_j)\bigr|
+\frac{c}{\sqrt{n}}
\int_0^{(2\sup_{g\in {\cal G}_{N+1}^{\prime}}P_n g^2)^{1/2}}
H_{d_{P_n,2}}^{1/2}({\cal G}_{N+1}^{\prime};u)du.
\label{ee3.3}
\end{equation}
Since we already proved (i)
it implies that on the event $E_{N+1}\bigcap {\cal J}$
\begin{eqnarray*}
&&
\inf_{g\in {\cal G}_{N+1}^{\prime}}
{\Bbb E}_{\eps}\bigl|n^{-1}\sum_{j=1}^n \eps_j g(X_j)\bigr|\leq
\inf_{g\in {\cal G}_{N+1}^{\prime}}
{\Bbb E}_{\eps}^{1/2}\bigl|n^{-1}\sum_{j=1}^n \eps_j g(X_j)\bigr|^2\leq
\frac{1}{\sqrt{n}}\inf_{g\in {\cal G}_{N+1}^{\prime}}\sqrt{P_n g^2}
\\
&&
\leq \frac{1}{\sqrt{n}}
\inf_{f\in {\cal F}_{N+1}}\sqrt{P_n \{f\leq \delta_{N,\frac{1}{2}}\}}\leq
\frac{1}{\sqrt{n}}
\inf_{f\in {\cal F}_{N+1}}\sqrt{P_n \{f\leq \delta\}}\leq
\sqrt{\frac{\eps}{n}}\leq \eps
\end{eqnarray*}
By the induction assumption, we also have on the event
$E_{N+1}\bigcap {\cal J}$
$$
\sup_{g\in {\cal G}_{N+1}^{\prime}}P_n g^2\leq \sup_{f\in {\cal F}_N}
P_n\{f\leq \delta_{N,\frac{1}{2}}\}\leq r_N.
$$
The bound for the Lipschitz norm of $\varphi_k^{\prime}$
gives the following bound on the distance
$$
d_{P_n,2}^2\Bigl(\varphi_{N+1}^{\prime}\circ f;
\varphi_{N+1}^{\prime}\circ g\Bigr)
=
n^{-1}\sum_{j=1}^n \Bigl|\varphi_{N+1}^{\prime} \circ f(X_j)-
\varphi_{N+1}^{\prime} \circ g(X_j) \Bigr|^2\leq
\Bigl(\frac{2}{\delta}\sqrt{\frac{r_{N}}{\eps}}\Bigr)^2 
d_{P_n,2}^2 (f,g).
$$
Therefore, on the event $E_{N+1}\bigcap {\cal J},$
we get quite similarly to (\ref{e3.6})
\begin{eqnarray}
&&
\frac{1}{\sqrt{n}}
\int_0^{(2\sup_{g\in {\cal G}_{N+1}^{\prime}}P_n g^2)^{1/2}}
H_{d_{P_n,2}}^{1/2}({\cal G}_{N+1}^{\prime};u)du
\leq
\frac{1}{\sqrt{n}}
\int_0^{(2r_N)^{1/2}}
H_{d_{P_n,2}}^{1/2}({\cal F};
\frac{\delta \sqrt{\eps} u}{2\sqrt{r_N}})du
\nonumber
\\
&&
\leq
c(\frac{r_N}{\eps})^{\alpha /4}
\frac{r_N^{1/2-\alpha/4}}{\sqrt{n}\delta^{\alpha/2}}
\leq c\sqrt{r_N \eps}.
\label{ee3.6}
\end{eqnarray}
We collect all bounds to see that on the event $E_{N+1}\bigcap {\cal J}$
\begin{equation}
\sup_{f\in {\cal F}_{N+1}}P_n\{f\leq \delta_{N+1}\}\leq
\frac{r_{N+1}}{2}+
c\sqrt{r_N \eps}.
\end{equation}
Therefore, it follows that with a proper choice of
constant $C>0$ in the recurrent relationship defining
the sequence $\{r_k\},$ we have on the event
$E_{N+1}\bigcap {\cal J}$
$$
\sup_{f\in {\cal F}_{N+1}}P_n\{f\leq \delta_{N+1}\}\leq
C\sqrt{r_N \eps}=r_{N+1},
$$
which proves the induction step for (ii)
and, therefore, the lemma
is proved.
\qed

To complete the proof of the theorem, we have to note that
the choice of
$N=[\log_2\log_2 \eps^{-1}]$
implies that $r_{N+1}\leq c\eps$ for some $c>0.$
The second inequality of the theorem can be proved
similarly with some minor modifications.
\qed

{\bf Proof of Theorem 5}.
Consider sequences
$\delta_j:=2^{-j\frac{2}{\gamma }},$
$$
\eps_j:=
\bigl(
\frac{1}{n\delta_j^{\alpha^{\prime}}}\bigr)^{\frac{2}{2+\alpha^{\prime}}},
\ j\geq 0,
$$
where $\alpha^{\prime}:=\frac{2\gamma}{2-\gamma}\geq \alpha .$
The first inequality of Theorem 6 implies
$$
{\Bbb P}\Bigl\{\exists j\geq 0\ \exists f\in {\cal F}\
P_n\{f\leq \delta_j\}\leq
\eps_j\
{\rm and}\
P\{f\leq \delta_j /2\}\geq  A^{\prime}\eps_j
\Bigr\}\leq
$$
\begin{equation}
\leq B^{\prime} \log_2 \log_2 n
\sum_{j\geq 0 }\exp\{-\frac{n^{\gamma/2}}{2}2^{2j}\}\leq
B \log_2 \log_2 n \exp\{-\frac{n^{\gamma/2}}{2}\}
\label{tristar}
\end{equation}
with some $B,B^{\prime},A^{\prime}>0$.
If for some $j\geq 1,$ we have
$$
\hat \delta_n(\gamma ;f)\in (\delta_j, \delta_{j-1}],
$$
then by definition of $\hat \delta_n(\gamma ;f)$
$$
P_n\{f\leq \delta_j\}\leq \eps_j.
$$
Suppose that for some $f\in {\cal F}$ the inequality 
$A^{-1}\hat\delta_n(\gamma;f)\leq \delta_n(\gamma;f)$
fails. Then, it follows
from the definition of $\delta_n(\gamma ;f)$ that
$$
P\{f\leq \delta_j /2\}\geq
P\Bigl\{f\leq \frac{\delta_{j-1}}{A}\Bigr\}
\geq
(\frac{1}{n\delta_{j-1}^{\alpha^{\prime}}})^{\frac{2}{2+\alpha^{\prime}}}
A^{\frac{2\alpha^{\prime}}{2+\alpha^{\prime}}}
\geq A^{\prime}\eps_j,
$$
where the last inequality holds for the proper
choice of a constant $A.$ Hence, (\ref{tristar})
guarantees the probability bound for the left side inequality
of the theorem. The right side inequality is proved similarly
utilizing the second inequality of Theorem 6.

\qed

\mdsk

\section{Convergence rates of empirical margin distributions}

As we defined in Section 2, ${\cal F}$ is a class of
measurable functions from $S$ into ${\Bbb R}.$
For $f\in {\cal F},$ let
$$
F_f(y):=P\{f\leq y\},\ F_{n,f}(y):=P_n\{f\leq y\},\ y\in {\Bbb R}.
$$
Let $L$ denote the L\'evy distance between the distribution functions in ${\Bbb R}:$
$$
L(F,G):=\inf \{\delta>0: F(t)\leq G(t+\delta)+\delta\ {\rm and}\ G(t)\leq F(t+\delta)+\delta ,\ {\rm for\ all
}\ t\in {\Bbb R}\}.
$$

In what follows, for a function $f$ from $S$ into ${\Bbb R}$ and $M>0,$
we denote $f_M$ the function that is equal to $f$ if $|f|\leq M,$
is equal to $M$ if $f>M$ and is equal to $-M$ if $f<-M.$ We set
$$
{\cal F}_M := \{f_M: f\in {\cal F}\}.
$$
As always, a function $F$ from $S$ into $[0,+\infty)$ is called
an envelope of ${\cal F}$ iff $|f(x)|\leq F(x)$ for all $f\in {\cal F}$
and all $x\in S.$

We write ${\cal F}\in GC(P)$ iff ${\cal F}$ is a Glivenko-Cantelli
class with respect to $P$
(i.e. $\|P_n-P\|_{\cal F}\to 0$ as $n\to\infty$ a.s.).
We write ${\cal F}\in BCLT(P)$ and say that ${\cal F}$
satisfies the Bounded Central Limit Theorem for $P$
iff
$$
{\Bbb E}\|P_n-P\|_{\cal F}=O(n^{-1/2}).
$$
In particular, this holds if ${\cal F}$ is a $P$-Donsker
class (see Dudley (1999), van der Vaart and Wellner (1996)
for precise definitions).

Our main goal in this section is to prove the following
results.

\begin{theorem}
Suppose that
\begin{equation}
\sup_{f\in {\cal F}}P\{|f|\geq M\}\to 0\ {\rm as}\ M\to\infty.
\label{tail1}
\end{equation}
Then, the following two statements are
equivalent:
$${\cal F}_M\in GC(P)\ {\rm for\ all\ } M>0\  \leqno (i)$$
and
$$
\sup_{f\in {\cal F}} L(F_{n,f},F_f)\to 0\ {\rm a.s.}\ {\rm as}\ n\goin .
\leqno (ii)
$$
\end{theorem}

\mdsk

\begin{theorem}
The following two statements are
equivalent:

(i) ${\cal F}\in GC(P)$

(ii) there exists a $P$-integrable envelope for the class 
${\cal F}^{(c)}=\{f-Pf : f\in{\cal F}\}$
and
$$
\sup_{f\in {\cal F}} L(F_{n,f},F_f)\to 0\ {\rm a.s.}\ {\rm as}\ n\goin .
$$
\end{theorem}

\mdsk
 
\begin{theorem}
Suppose that the class ${\cal F}$ is uniformly bounded.
If ${\cal F}\in BCLT(P),$ then
$$
\sup_{f\in {\cal F}}L(F_{n,f},F_f)=O_P(n^{-1/4})\ {\rm as}\ n\goin .
$$
Moreover, if for some $\alpha \in (0,2)$ and for some $D>0$
\begin{equation}
H_{d_{P_n,2}}({\cal F};u)\leq Du^{-\alpha},\ u>0\ {\rm a.s.},
\label{entro}
\end{equation}
then
$$
\sup_{f\in {\cal F}}L(F_{n,f},F_f)=O(n^{-\frac{1}{2+\alpha}})\ {\rm as}\ n\goin \ {\rm a.s.}
$$
\end{theorem}

The following theorem gives the bound that plays an important
role in the proofs.

\begin{theorem}
Let $M>0$ and let ${\cal F}$ be a class of measurable functions from
$S$ into $[-M,M].$
For all $t>0,$
$$
{\Bbb P}\Bigl\{\sup_{f\in {\cal F}}L(F_{n,f},F_f)\geq
2\Bigl({\Bbb E}\|n^{-1}\sum_{i=1}^n \eps_i\delta_{X_i}\|_{\cal
F}+
{{M}\over {\sqrt {n}}}\Bigr)^{1/2}+\frac{t}{\sqrt{n}}
\Bigr\}\leq \exp \{-2t^2\}.
$$
\end{theorem}

{\bf Proof}. Let $\delta >0.$
Let $\varphi(x)$ be equal to $1$ for $x\leq 0,$
$0$ for $x\geq 1$ and linear in between.
One can get the following bounds:
$$
F_f(y)=P\{f\leq y\} \leq P\varphi ({{f-y}\over {\delta}})\leq
P_n \varphi ({{f-y}\over {\delta}})+
\|P_n-P\|_{\tilde {\cal G}_{\delta}}
$$
$$
\leq F_{n,f}(y+\delta)
+\|P_n-P\|_{\tilde {\cal G}_{\delta}}
$$
and
$$
F_{n,f}(y)=P_n\{f\leq y\} \leq P_n\varphi ({{f-y}\over {\delta}})\leq
P \varphi ({{f-y}\over {\delta}})+\|P_n-P\|_{\tilde {\cal G}_{\delta}}
$$
$$
\leq
F_{f}(y+\delta)+\|P_n-P\|_{\tilde {\cal G}_{\delta}},
$$
where
$$
\tilde {\cal G}_{\delta}:=
\Bigl\{\varphi \circ ({{f-y}\over \delta}) - 1:
f\in {\cal F}, y\in [-M,M]\Bigr\}.
$$
Similarly to the proof of Theorem 1
we get that with probability at least $1-2e^{-2t^2}$
\begin{equation}
\|P_n-P\|_{\tilde {\cal G}_{\delta}}
\leq
\frac{4}
{\delta}\Bigl[{\Bbb E}\|n^{-1}\sum_{i=1}^n \eps_i\delta_{X_i}\|_{\cal F}+
M n^{-1/2}\Bigr]+\frac{t}{\sqrt{n}}. \label{e4.2}
\end{equation}
Setting
$$
\delta :=
2\Bigl({\Bbb E}\|n^{-1}\sum_{i=1}^n \eps_i\delta_{X_i}\|_{\cal
F}+ Mn^{-1/2}\Bigr)^{1/2},
$$
we get that
with probablity at least $1-\exp \{-2t^2\}$
$$
\sup_{f\in {\cal F}} L(F_{n,f},F_f) \leq
2\Bigl({\Bbb E}\|n^{-1}\sum_{i=1}^n \eps_i\delta_{X_i}\|_{\cal
F}+ M n^{-1/2}\Bigr)^{1/2}+\frac{t}{\sqrt{n}},
$$
which completes the proof.
\qed

{\bf Proof of Theorem 7}.
First we prove that (i) implies (ii).
Since ${\cal F}_M\in GC(P),$ we have
$${\Bbb E}\|P_n-P\|_{{\cal F}_M}\to 0\ {\rm as}\ n\goin,$$
which, by symmetrization inequality, implies
$${\Bbb E}\|n^{-1}\sum_{i=1}^n \eps_i \delta_{X_i}\|_{{\cal F}_M}
\to 0\ {\rm as}\ n\goin .
$$
Plugging in the bound of Theorem 10 $t=\log n$ and using Borel-Cantelli
Lemma proves that for all $M>0$
$$
\sup_{f\in {\cal F}}L(F_{n,f_M},F_{f_M})=
\sup_{f\in {\cal F}_M}L(F_{n,f},F_f)\to 0\ {\rm as}\ n\to\infty\ {\rm a.s.}
$$
The following bounds easily follow from the definition of
L\'evy distance:
$$
\sup_{f\in {\cal F}}L(F_f,F_{f_M})\leq
\sup_{f\in {\cal F}}P\{|f|\geq M\}
$$
and
$$
\sup_{f\in {\cal F}}L(F_{n,f},F_{n,f_M})\leq
\sup_{f\in {\cal F}}P_n\{|f|\geq M\}.
$$
By condition (\ref{tail1}) of the theorem,
$$
\sup_{f\in {\cal F}}L(F_f,F_{f_M})\to 0\ {\rm as}\ M\to \infty .
$$
To prove that also
$$
\lim_{M\goin}\limsup_{n\goin}\sup_{f\in {\cal F}}
L(F_{n,f},F_{n,f_M})=0\ {\rm a.s.},
$$
it is enough to show that
\begin{equation}
\lim_{M\goin}\limsup_{n\goin}\sup_{f\in {\cal F}}
P_n\{|f|\geq M\}=0\ {\rm a.s.}
\label{whatleft}
\end{equation}
To this end, consider the function $\varphi$ from
${\Bbb R}$ into $[0,1]$ that is equal to $0$ for
$|u|\leq M-1,$ is equal to $1$ for $|u|>M$ and is
linear in between. We have
\begin{eqnarray}
&&
\sup_{f\in {\cal F}}P_n\{|f|\geq M\}=
\sup_{f\in {\cal F}_M}P_n\{|f|\geq M\}\leq
\sup_{f\in {\cal F}_M}P_n \varphi (|f|)
\nonumber
\\
&&
\leq
\sup_{f\in {\cal F}_M}P\varphi (|f|)+\|P_n-P\|_{\cal G}\leq
\sup_{f\in {\cal F}_M}P\{|f|\geq M-1\}+\|P_n-P\|_{\cal G},
\label{whatis}
\end{eqnarray}
where
$$
{\cal G}:=\Bigl\{\varphi \circ f: f\in {\cal F}_M\Bigr\}.
$$
Since $\varphi$ satisfies the Lipschitz condition with
constant $1,$ the argument based on symmetrization inequality
and comparison inequalities
(see the proofs above) allows one to show that the
condition (i) implies that
$$
{\Bbb E}\|P_n-P\|_{\cal G} \to 0\ {\rm as}\ n\goin .
$$
Then, the standard use of concentration inequality implies that
$$
\|P_n-P\|_{\cal G} \to 0\ {\rm as}\ n\goin \ {\rm a.s.}
$$
Therefore, (\ref{whatleft}) immediately follows from
condition (\ref{tail1}) and (\ref{whatis}).
Now, the triangle inequality for the L\'evy distance
allows one easily to complete the proof of (ii).

To prove that (ii) implies (i), we use the following bound
$$
|\int_{-M}^{M} td(F-G)(t)|\leq c L(F,G),
$$
which holds with some constant $c=c(M)$ for any two distribution
functions on $[-M,M].$ The bound implies that
\begin{equation}
\|P_n-P\|_{{\cal F}_M}=\sup_{f\in {\cal F}_M}|P_n f-Pf|=
\sup_{f\in {\cal F}_M}|\int_{-M}^{M}td(F_{n,f}-F_f)(t)|
\leq c\sup_{f\in {\cal F}_M}L(F_{n,f};F_f).
\label{levy1}
\end{equation}
Since for all $M>0$ and for all $f\in {\cal F}$ it is easily
proved that
\begin{equation}
L(F_{n,f_M},F_{f_M})\leq L(F_{n,f},F_f),
\label{LM}
\end{equation}
the bound (\ref{levy1}) and condition (ii) imply (i),
which completes the proof of the second statement.
\qed

\mdsk

{\bf Proof of Theorem 8}. 
Since centering does not change L\'evy distance and
does not change Glivenko-Cantelli property we can start
by assuming that $\cal F$ is centered, i.e. 
${\cal F}={\cal F}^{(c)}.$
To prove that (i) implies (ii),
note first of all that
the condition ${\cal F}\in GC(P)$ yields that 
${\cal F}= {\cal F}^{(c)}$
has a $P$-integrable envelope (see van der Vaart and Wellner (1996),
p. 125). Also, the existence of a $P$-integrable
envelope implies (\ref{tail1}).
Finally, if ${\cal F}\in GC(P),$ then for all $M>0$
${\cal F}_M\in GC(P)$
[To prove this claim note that $f_M=\varphi_M \circ f,$
where $\varphi_M$ is the function from ${\Bbb R}$ into
$[-M,M]$ that is equal to $u$ for $|u|\leq M,$ is equal
to $M$ for $u> M$ and is equal to $-M$ for $u<-M.$
The function $\varphi_M$ is Lipschitz with constant $1$
which allows to prove the claim by the argument based
on the comparison inequality and used many times above].
We can use Theorem 7 to conclude that (i) implies (ii).
On the other hand, if (ii) holds then by the inequality
(\ref{LM}) we get that
$$
\sup_{f\in {\cal F}_M}L(F_{n,f},F_f)\to 0\ {\rm as}\ n\goin\ {\rm a.s.}
$$
As we pointed out above (\ref{tail1}) holds, so, by Theorem 7,
we have ${\cal F}_M\in GC(P)$ for all $M>0.$ The integrability
of the envelope of the class ${\cal F}$ allows us to conclude
the proof of (i) by a standard truncation argument.
\qed

\mdsk

{\bf Proof of Theorem 9}. Since ${\cal F}$ is uniformly bounded,
we can choose $M>0$ such that ${\cal F}_M={\cal F}.$
To prove the first statement note that ${\cal F}\in BCLT(P)$ means that
$$
{\Bbb E} \|P_n-P\|_{\cal F}=O(n^{-1/2}).
$$
which implies
$$
{\Bbb E}\|n^{-1}\sum_{i=1}^n \eps_i \delta_{X_i}\|_{\cal F}=O(n^{-1/2}).
$$
Thus, the bound of Theorem 10 implies that with some constant $C>0$
$$
{\Bbb P}\{\sup_{f\in {\cal F}}L(F_{n,f},F_f)\geq
\bigl({C\over {\sqrt{n}}}+{{4M}\over {\sqrt{n}}}\bigr)^{1/2}
+\frac{t}{\sqrt{n}}\}\leq \exp\{-2t^2\}.
$$
It follows that
$$
\lim_{u\goin}\limsup_{n\goin}{\Bbb P}\{n^{1/4}\sup_{f\in {\cal F}}L(F_{n,f},F_f)\geq
u\}=0.
$$

To prove the second statement, we follow the proof
of Theorem 10. We use Rademacher symmetrization
inequality to get the bound
$$
{\Bbb E}\|P_n-P\|_{\tilde {\cal G}_{\delta}}\leq
2{\Bbb E}\hat R_n(\tilde {\cal G}_{\delta})
$$
and then use the entropy inequalities for subgaussian
processes (see \cite{Well}, Corollary 2.2.8)
to show that
\begin{eqnarray*}
&&
{\Bbb E}_{\eps}\hat R_n(\tilde {\cal G}_{\delta})\leq
\inf_{g\in \tilde{\cal G}_{\delta}}
{\Bbb E}_{\eps}\Bigl|n^{-1}\sum_{j=1}^n\eps_j g(X_j)\Bigr|+
\frac{c}{\sqrt{n}}
\int_0^{\sqrt{2\sup_{g\in \tilde{\cal G}_{\delta}}P_ng^2}}
H_{d_{P_n,2}}^{1/2}\Bigl(\tilde{\cal G}_{\delta};u\Bigr)du
\\
&&
\leq
\frac{1}{\sqrt{n}}+
\frac{c}{\sqrt{n}}
\int_0^{\sqrt{2}}
H_{d_{P_n,2}}^{1/2}\Bigl(\tilde{\cal G}_{\delta};u\Bigr)du.
\end{eqnarray*}
To bound the random entropy $H_{d_{P_n,2}},$ we use
the Lipschitz condition for the function $\varphi.$
It yields (via a standard argument based on constructing
minimal covering of the class ${\cal F}$ with respect to the metric
$d_{P_n,2}$ and of the interval $[-M,M]$ with respect to the usual distance
in real line and "combining" the coverings properly) the following
bound:
$$
H_{d_{P_n,2}}\Bigl(\tilde{\cal G}_{\delta};u\Bigr)\leq
H_{d_{P_n,2}}\Bigl({\cal F};\delta u/2\Bigr)+\log{\frac{4M}{u\delta}}.
$$
Therefore, we get (with a proper constant $c>0$)
$$
{\Bbb E}_{\eps}\hat R_n(\tilde {\cal G}_{\delta})\leq
\frac{c}{\sqrt{n}}\Bigr[
\int_0^{\sqrt{2}}
H_{d_{P_n,2}}^{1/2}\Bigl({\cal F};\delta u\Bigr)du+
\sqrt{\log{\frac{4M}{\delta}}}+1\Bigl],
$$
which, under the condition (\ref{entro}), is bounded
from above by
$\frac{c}{\sqrt{n}\delta^{\alpha/2}}.$
Thus, we proved the bound
$$
{\Bbb E}\|P_n-P\|_{\tilde {\cal G}_{\delta}}
\leq
\frac{c}{\sqrt{n}\delta^{\alpha/2}}.
$$
Arguing now the same way as in the proof of Theorem 10,
we can show that with probability at least $1-\exp\{-2t^2\},$
$$
\sup_{f\in {\cal F}}L(F_{n,f},F_f)\leq \delta \bigvee
\frac{c}{\sqrt{n}\delta^{\alpha/2}}+\frac{t}{\sqrt{n}}.
$$
Plugging in the last inequality
$$
\delta := \frac{c}{n^{\frac{1}{2+\alpha}}},
$$
we get
$$
{\Bbb P}\Bigl\{\sup_{f\in {\cal F}}L(F_{n,f},F_f)\geq
\frac{c}{n^{\frac{1}{2+\alpha}}}+\frac{t}{\sqrt{n}}\Bigr\}\leq
\exp\{-2t^2\}.
$$
By choosing $t:=\log {n}$ and using Borel-Cantelli Lemma,
we complete the proof of the second statement.

\qed

\mdsk


{\bf Remark}. It's interesting to mention that the condition
${\cal F}\in GC(P)$ \it does not \rm imply that
$$
\sup_{f\in \cal F} \sup_{t\in {\Bbb R}}|F_{n,f}(t)-F_f(t)|\to 0
$$
with probability 1,
which is equivalent to saying that the class of sets
$\{I(f\leq t) : f\in{\cal F}, t\in {\Bbb R} \}$
is $GC(P)$.
As an example, consider the case when
$S$ is a unit ball in an infinite-dimensional separable
Banach space. Let ${\cal F}$ be the restriction of the unit
ball in the dual space on $S$. For i.i.d. random variables $\{X_n\}$
in $S,$ we have, by the LLN in separable Banach spaces,
$$
\|P_n-P\|_{\cal F}:=\Bigl\|n^{-1}\sum_{j=1}^n (X_j-{\Bbb E}X)\Bigr\|
\to 0\ {a.s.},
$$
so ${\cal F}\in GC(P).$ On the other hand, there exists an example
of a distribution $P$ such that ${\cal H}\not\in GC(P),$
where ${\cal H}$ is the class of all halfspaces
(see Sazonov (1963) and also Tops\o e, Dudley and Hoffmann-J\o rgensen (1976)).
Hence,
$$
\sup_{f\in {\cal F}}\sup_{t\in {\Bbb R}}|F_{n,f}(t)-F_{f}(t)|=
\|P_n-P\|_{\cal H}
$$
does not converge to $0$ a.s.

In the next proposition, we are again considering the class ${\cal F}$ used already in
Proposition 1 and the sequence of observations $\{X_n\}$ defined by
$$
X_n=\left\{
\eps_k^n (2\log(k+1))^{-\frac{1}{2}-\beta}
\right\}_{k\geq 1},\ n\geq 1,
$$
where
$\beta:=\frac{1}{\alpha}-\frac{1}{2},$ $\alpha \in (0,2]$
and
$\eps_k^n$
are i.i.d. Rademacher random variables.
The proposition shows the optimality of the rates
of convergence obtained in Theorem 9.

\begin{proposition}
Consider the sequence
$\delta_n$
such that
$$
\sup_{f\in \cal F} L(F_{n,f},F_f)=O_{P}(\delta_n).
$$
Then
$$
\delta_n \geq cn^{-\frac{1}{2+\alpha}}
$$
(when
$\alpha =2,$ we have $\delta_n \geq cn^{-1/4}$).
On the other hand, for $\alpha \in (0,2),$ we have
$$
H_{d_{P_n,2}}({\cal F};u)\leq D u^{-\alpha},\ u>0
$$
and
$$
\sup_{f\in \cal F} L(F_{n,f},F_f)=
O(n^{-\frac{1}{2+\alpha}})\ {\rm a.s.};
$$
for $\alpha =2$ we have ${\cal F}\in BCLT(P)$ and
$$
\sup_{f\in \cal F} L(F_{n,f},F_f)=
O_P(n^{-\frac{1}{4}}).
$$
\end{proposition}

{\bf Proof}.
We can assume without loss of generality
that with probability more than $1/2$
for all
$k\geq 1,$ $y\in [-1,1]$ and $n$ large enough we have
\begin{equation}
P(f_k\leq y)\leq P_n(f_k\leq y+\delta)+\delta.\label{star}
\end{equation}
If we take
$y=0$ and consider only such $k$ that
satisfy the inequality
$(2\log(k+1))^{\beta+1/2}<\delta^{-1}$ then (\ref{star})
becomes equivalent to
$$
1/2\leq n^{-1}\sum_{i\leq n}I(\eps_k^i = -1)+\delta.
$$
Inequality
$(2\log(k+1))^{\beta+1/2}<\delta^{-1}$
holds for
$k\leq \psi_1(\delta)=1/2\exp(\delta^{-\frac{2}{1+2\beta}}/2).$
Therefore, for large $n$
\begin{eqnarray}
1/2
&\leq &
{\Bbb P}\Biggl\{
\bigcap_{k\leq \psi_1(\delta)}\Bigl\{
1/2\leq n^{-1}\sum_{i\leq n}I(\eps_k^i = -1)+\delta
\Bigr\}
\Biggr\}
\nonumber
\\
&=&
{\Bbb P}\Biggl\{
1/2\leq n^{-1}\sum_{i\leq n}I(\eps_1^i = -1)+\delta
\Biggr\}^{\psi_1(\delta)}
\leq
\Biggl(1-{n\choose k_0}2^{-n}\Biggr)^{\psi_1(\delta)},
\label{twostars}
\end{eqnarray}
where
$k_0=[n/2-\delta n]-1.$
Using (\ref{Lower}), we get
$$
2^{-\frac{1}{\psi_1(\delta)}}\leq
1-cn^{-\frac{1}{2}}\exp(-4n\delta^2).
$$
Taking logarithm of both sides and taking into account
that $\log(1-x)\leq -x$ we get
(recall that $\psi_1(\delta)=1/2\exp(\delta^{-\frac{2}{1+2\beta}}/2)$)
$$
\exp(-2^{-1}\delta^{-\frac{2}{1+2\beta}})\geq
cn^{-\frac{1}{2}}\exp(-4n\delta^2).
$$
Therefore,
$$
1/(2\delta^{2/(1+2\beta)})\leq 4n\delta^2+c\log n
$$
and
$$
1/2\leq 4n\delta^{4(1+\beta)/(1+2\beta)}+
c\delta^{2/(1+2\beta)}\log n.
$$
This finally implies that
$$
\delta\geq cn^{-\frac{1+2\beta}{4(1+\beta)}}=cn^{-\frac{1}{2+\alpha}}.
$$
The second statement follows from Theorem 9.
To check condition
(\ref{entro}), note that in this case,
as soon as $2\log N\geq (u/2)^{-\alpha},$ we
have $|f_k (X_n)|\leq u/2$ for all $k\geq N$
and $n\geq 1.$ Hence,
$$
d_{P_n,2}(f_k, f_N)\leq u,\ k\geq N
$$
and we have
$$
H_{d_{P_n,2}}({\cal F};u)\leq \log N,
$$
which implies (\ref{entro}).
For $\alpha=2,$ we also have ${\cal F}\in BCLT(P)$
(see Ledoux and Talagrand (1991), pp. 276--277).
Theorem 9 allows one to complete the proof.
\qed

\mdsk

\section{Bounding the generalization error of convex combinations of classifiers}

In this and in the next section we consider applications of
the bounds of Section 2 to various learning (classification)
problems.
We start with an application of the inequalities of
Section 2 to bounding the generalization error in general
multiclass problems. Namely, we assume that the labels take
values in a finite set ${\cal Y}$ with ${\rm card}({\cal Y})=M.$
Consider a class $\tilde {\cal F}$ of functions from ${\tilde S}:=S\times {\cal Y}$
into ${\Bbb R}.$ A function $f\in \tilde {\cal F}$ predicts a label $y\in {\cal Y}$
for an example $x\in S$ iff
$$
f(x,y)>\max_{y^{\prime}\neq y}f(x,y^{\prime}).
$$
The margin of a labeled example $(x,y)$ is defined as
$$
m_f(x,y):=f(x,y)-\max_{y^{\prime}\neq y}f(x,y^{\prime}),
$$
so $f$ misclassifies the labeled example $(x,y)$ iff $m_f(x,y)\leq 0.$
Let
$$
{\cal F}:= \{f(\cdot, y):y\in {\cal Y}, f\in \tilde{\cal F}\}.
$$
The proof of the next result is based on the application of Theorem 2.

\begin{theorem}
For all $t >0,$
$$
{\Bbb P}\Bigl\{\exists f\in \tilde {\cal F}:
P\{m_f\leq 0\} > \inf_{\delta \in (0,1]}\Bigl[P_n\{m_f\leq \delta\}
+{{8M(2M-1)}\over {\delta}}R_n({\cal F})
$$
$$
+\Bigl(\frac{\log\log_2 (2\delta^{-1})}{n}\Bigr)^{1/2}
\Bigr]
+{t\over {\sqrt {n}}}
\Bigr\}
\leq 2\exp \{-2t^2\}.
$$
\end{theorem}

To prove the theorem, we use the following lemma.

For a class of functions ${\cal H},$ we will denote by
$$
{\cal H}^{(l)}=\{
\max (h_1,\ldots,h_l): h_1,\ldots,h_l\in {\cal H}
\}.
$$

\begin{lemma} The following bound holds:
$$
{\Bbb E}\|\sum_{i=1}^n \eps_i\delta_{X_i}\|_{{\cal H}^{(l)}}\leq
2l{\Bbb E}\|\sum_{i=1}^n \eps_i\delta_{X_i}\|_{\cal H}.
$$
\end{lemma}

{\bf Proof}.
Let $x^{+}:=x\vee 0.$
Obviously $x\mapsto x^+$ is a nondecreasing convex function
such that $(a+b)^+ \leq a^+ +b^+ .$
We will first prove that
\begin{equation}
{\Bbb E}(\sup_{{\cal H}^{(l)}}
\sum_{i=1}^n \eps_i h(X_i))^+
\leq
l{\Bbb E}(\sup_{\cal H}
\sum_{i=1}^n \eps_i h(X_i))^+ .
\label{Medstep}
\end{equation}
Let us consider classes of functions
${\cal F}_1,$  ${\cal F}_2$ and
$$
{\cal F}=\{
\max (f_1,f_2): f_1\in {\cal F}_1,
f_2\in {\cal F}_2
\}.
$$
Since
$$
\max (f_1, f_2)=\frac{1}{2}\left(
(f_1+f_2)+|f_1-f_2|
\right),
$$
we have
\begin{eqnarray*}
&&
{\Bbb E}\Bigl(\sup_{{\cal F}}
\sum_{i=1}^n \eps_i f(X_i)\Bigr)^+
\\
&&
\leq
{\Bbb E}\Bigl(\sup_{{\cal F}_1,{\cal F}_2}
\sum_{i=1}^n \eps_i\frac{1}{2}(f_1(X_i)+f_2(X_i))
+\sup_{{\cal F}_1,{\cal F}_2}
\sum_{i=1}^n \eps_i\frac{1}{2}|f_1(X_i)-f_2(X_i)|
\Bigr)^+
\\
&&
\leq
\frac{1}{2}
{\Bbb E}\Bigl(\sup_{{\cal F}_1,{\cal F}_2}
\sum_{i=1}^n \eps_i(f_1(X_i)+f_2(X_i))\Bigr)^+
+\frac{1}{2}
{\Bbb E}\Bigl(\sup_{{\cal F}_1,{\cal F}_2}
\sum_{i=1}^n \eps_i|f_1(X_i)-f_2(X_i)|\Bigr)^+
\\
&&
\leq
\frac{1}{2}
{\Bbb E}\Bigl(\sup_{{\cal F}_1}
\sum_{i=1}^n \eps_i f_1(X_i)\Bigr)^+
+
\frac{1}{2}
{\Bbb E}\Bigl(\sup_{{\cal F}_2}
\sum_{i=1}^n \eps_i f_2(X_i)\Bigr)^+ 
\\
&&
+\frac{1}{2}
{\Bbb E}\Bigl(\sup_{{\cal F}_1,{\cal F}_2}
\sum_{i=1}^n \eps_i|f_1(X_i)-f_2(X_i)|\Bigr)^+ .
\end{eqnarray*}
The proof of Theorem 4.12 in \cite{LT}
contains the following statement. If $T$ is a bounded subset
of ${\Bbb R}^n,$ functions $\varphi_i,\,i=1,\ldots,n$
are contractions such that $\varphi_i(0)=0$
and a function $G:{\Bbb R}\to {\Bbb R}$ is convex and
nondecreasing then
$$
{\Bbb E}G\Bigl(
\sup_{t\in T}\sum_{i=1}^{n}\eps_i \varphi_i(t_i)
\Bigr)\leq
{\Bbb E}G\Bigl(
\sup_{t\in T}\sum_{i=1}^{n}\eps_i t_i
\Bigr).
$$
If we take $G(x)=x^+,$ $\varphi_i(x)=|x|$
and $T=\{(f_1(X_i)-f_2(X_i))_{i=1}^n :
f_1\in {\cal F}_1, f_2\in {\cal F}_2\}$ we get
(first conditionally on $(X_i)_{i=1}^n$ and then taking
expectations)
\begin{eqnarray*}
&&
{\Bbb E}\Bigl(\sup_{{\cal F}_1,{\cal F}_2}
\sum_{i=1}^n \eps_i|f_1(X_i)-f_2(X_i)|\Bigr)^+
\leq
{\Bbb E}\Bigl(\sup_{{\cal F}_1,{\cal F}_2}
\sum_{i=1}^n \eps_i(f_1(X_i)-f_2(X_i))\Bigr)^+
\\
&&
\leq
{\Bbb E}\Bigl(\sup_{{\cal F}_1}
\sum_{i=1}^n \eps_i f_1(X_i)\Bigr)^+
+{\Bbb E}\Bigl(\sup_{{\cal F}_2}
\sum_{i=1}^n \eps_i f_2(X_i)\Bigr)^+ ,
\end{eqnarray*}
where in the last inequality we used the fact
that the sequence
$(-\eps_i)_{i=1}^n$
is equal in distribution to
$(\eps_i)_{i=1}^n.$
Combining the bounds gives
$$
{\Bbb E}\Bigl(\sup_{{\cal F}}
\sum_{i=1}^n \eps_i f(X_i)\Bigr)^+
\leq
{\Bbb E}\Bigl(\sup_{{\cal F}_1}
\sum_{i=1}^n \eps_i f_1(X_i)\Bigr)^+
+{\Bbb E}\Bigl(\sup_{{\cal F}_2}
\sum_{i=1}^n \eps_i f_2(X_i)\Bigr)^+ .
$$
Now by induction we easily get (\ref{Medstep}).
Finally, again using the fact that $(-\eps_i)_{i=1}^n$
is equal in distribution to
$(\eps_i)_{i=1}^n,$ we conclude the proof:
\begin{eqnarray*}
&&
{\Bbb E}\|\sum_{i=1}^n \eps_i\delta_{X_i}\|_{{\cal H}^{(l)}}
\leq
{\Bbb E}\Bigl(\sup_{{\cal H}^{(l)}}
\sum_{i=1}^n \eps_i h(X_i)\Bigr)^+
+{\Bbb E}\Bigl(-\sup_{{\cal H}^{(l)}}
\sum_{i=1}^n \eps_i h(X_i)\Bigr)^+ 
\\
&&
=2{\Bbb E}\Bigl(\sup_{{\cal H}^{(l)}}
\sum_{i=1}^n \eps_i h(X_i)\Bigr)^+
\leq
2l{\Bbb E}\Bigl(\sup_{{\cal H}}
\sum_{i=1}^n \eps_i h(X_i)\Bigr)^+
\leq
2l{\Bbb E}\|\sum_{i=1}^n \eps_i\delta_{X_i}\|_{{\cal H}}.
\end{eqnarray*}
\qed

{\bf Proof of Theorem 11.} We have the following bounds:
\begin{eqnarray*}
&&
{\Bbb E}\sup_{f\in \tilde
{\cal F}}\Bigl|n^{-1}\sum_{j=1}^n \eps_j m_f(X_j,Y_j)\Bigr|=
{\Bbb E}\sup_{f\in \tilde {\cal F}}
\Bigl|n^{-1}\sum_{j=1}^n \eps_j
\sum_{y\in {\cal Y}} m_f(X_j,y)I_{\{Y_j=y\}}\Bigr|
\\
&&
\leq \sum_{y\in {\cal Y}}{\Bbb E}\sup_{f\in \tilde {\cal F}}
\Bigl|n^{-1}\sum_{j=1}^n \eps_j m_f(X_j,y)I_{\{Y_j=y\}}\Bigr|
\\
&&
\leq {1\over 2}\sum_{y\in {\cal Y}}{\Bbb E}\sup_{f\in \tilde {\cal F}}
\Bigl|n^{-1}\sum_{j=1}^n \eps_j m_f(X_j,y)(2I_{\{Y_j=y\}}-1)\Bigr| +
{1\over 2}\sum_{y\in {\cal Y}}{\Bbb E}\sup_{f\in \tilde {\cal F}}
\Bigl|n^{-1}\sum_{j=1}^n \eps_j m_f(X_j,y)\Bigr|.
\end{eqnarray*}
Denote $\sigma_j(y):=2I_{\{Y_j=y\}}-1.$ Given $\{(X_j,Y_j):1\leq j\leq n\},$
the random variables $\{\eps_j\sigma_j(y):1\leq j\leq n\}$ are
i.i.d. Rademacher.
Hence, we have
\begin{eqnarray*}
&&
{\Bbb E}\sup_{f\in \tilde {\cal F}}
\Bigl|n^{-1}\sum_{j=1}^n \eps_j m_f(X_j,y)(2I_{\{Y_j=y\}}-1)\Bigr| =
{\Bbb E}\sup_{f\in \tilde {\cal F}}
\Bigl|n^{-1}\sum_{j=1}^n \eps_j \sigma_j(y) m_f(X_j,y)\Bigr|
\\
&&
={\Bbb E}{\Bbb E}_\eps\sup_{f\in \tilde {\cal F}}
\Bigl|n^{-1}\sum_{j=1}^n \eps_j \sigma_j(y) m_f(X_j,y)\Bigr| =
{\Bbb E}{\Bbb E}_\eps\sup_{f\in \tilde {\cal F}}
\Bigl|n^{-1}\sum_{j=1}^n \eps_j  m_f(X_j,y)\Bigr| 
\\
&&
={\Bbb E}\sup_{f\in \tilde {\cal F}}
\Bigl|n^{-1}\sum_{j=1}^n \eps_j  m_f(X_j,y)\Bigr| .
\end{eqnarray*}
Therefore, we have
$$
{\Bbb E}\sup_{f\in \tilde {\cal F}}
\Bigl|n^{-1}\sum_{j=1}^n \eps_j m_f(X_j,Y_j)\Bigr|\leq
\sum_{y\in {\cal Y}}
{\Bbb E}\sup_{f\in \tilde {\cal F}}
\Bigl|n^{-1}\sum_{j=1}^n \eps_j m_f(X_j,y)\Bigr|.
$$
Next, using Lemma 2, we get for all $y\in {\cal Y}$
$$
{\Bbb E}\sup_{f\in \tilde {\cal F}}
\Bigl|n^{-1}\sum_{j=1}^n \eps_j m_f(X_j,y)\Bigr|\leq
{\Bbb E}\sup_{f\in \tilde {\cal F}}
\Bigl|n^{-1}\sum_{j=1}^n \eps_j f(X_j,y)\Bigr|+
{\Bbb E}\sup_{f\in \tilde {\cal F}}
\Bigl|n^{-1}\sum_{j=1}^n
\eps_j \max_{y^{\prime}\neq y}f(X_j,y^{\prime})\Bigr|
$$
$$
\leq
{\Bbb E}\sup_{f\in {\cal F}}
\Bigl|n^{-1}\sum_{j=1}^n \eps_j f(X_j)\Bigr|+
{\Bbb E}\sup_{f\in {\cal F}^{(M-1)}}
\Bigl|n^{-1}\sum_{j=1}^n \eps_j f(X_j)\Bigr|
$$
$$
\leq
(2M-1) {\Bbb E}\sup_{f\in {\cal F}}
\Bigl|n^{-1}\sum_{j=1}^n \eps_j f(X_j)\Bigr|.
$$
This implies
\begin{eqnarray*}
&&
{\Bbb E}\sup_{f\in \tilde {\cal F}}
\Bigl|n^{-1}\sum_{j=1}^n \eps_j m_f(X_j,Y_j)\Bigr|
\leq 
\sum_{y\in {\cal Y}}
(2M-1){\Bbb E}\sup_{f\in {\cal F}}
\Bigl|n^{-1}\sum_{j=1}^n \eps_j f(X_j)\Bigr|
\\
&&
=
M(2M-1) {\Bbb E}\sup_{f\in {\cal F}}
\Bigl|n^{-1}\sum_{j=1}^n \eps_j f(X_j)\Bigr|,
\end{eqnarray*}
and the result follows from Theorem 2
(one can use in this theorem the continuous function $\varphi$ that is
equal to $1$ on $(-\infty, 0],$ is equal to $0$ on $[1,+\infty)$
and is linear in between).
\qed

In the rest of the paper, we assume that the set of labels is
$\{-1,1\},$ so that
$\tilde S:=S\times \{-1,1\}$ and
${\tilde {\cal F}}:=\{\tilde f: f\in {\cal F}\},$ where $\tilde f(x,y):=yf(x).$
$P$ will denote the distribution of $(X,Y),$ $P_n$ the empirical distribution
based on the observations $((X_1,Y_1),\dots, (X_n,Y_n)).$
Clearly, we have
$$
R_n({\tilde {\cal F}})=
{\Bbb E}\sup_{f\in {\cal F}}\bigl|n^{-1}\sum_{i=1}^n \eps_i Y_i f(X_i)\bigr|=
{\Bbb E}
{\Bbb E}_{\eps}
\sup_{f\in {\cal F}}\bigl|n^{-1}\sum_{i=1}^n {\tilde \eps_i}  f(X_i)\bigr|,
$$
where $\tilde \eps_i := Y_i \eps_i.$
Since, for given $\{(X_i,Y_i)\},$ $\{\tilde \eps_i\}$ and $\{\eps_i\}$
have the same distribution, we get
$$
{\Bbb E}_{\eps}
\sup_{f\in {\cal F}}\bigl|n^{-1}\sum_{i=1}^n {\tilde \eps_i}  f(X_i)\bigr|=
{\Bbb E}_{\eps}\sup_{f\in {\cal F}}\bigl|n^{-1}\sum_{i=1}^n
\eps_i  f(X_i)\bigr|,
$$
which immediately implies $R_n({\tilde {\cal F}})=R_n({\cal F}).$

The results of Section 2 now give some useful bounds for boosting and other methods of
combining the classifiers. Namely, we get in this case the following
theorem (compare with the recent result of
Schapire, Freund, Bartlett and Lee (1998)).

Given a class ${\cal H}$ of measurable functions from $S$ into
${\Bbb R},$ we denote ${\rm conv}({\cal H})$ the closed convex
hull of ${\cal H},$ i.e. ${\rm conv }({\cal H})$ consists of
all functions on $S$ that are pointwise limits of convex
combinations of functions from ${\cal H}:$
\begin{eqnarray*}
&&
{\rm conv}({\cal H}):=\Bigl\{
f:\ \forall x\in S\
f(x)=\lim f_N(x),\ 
f_N=\sum_{j=1}^{N} 
w_j^N h_j^N,
\\
&&
w_j^N\geq 0,\ \sum_{j=1}^{N} w_j^N=1,\ h_j^N\in {\cal H},\ N\geq 1
\Bigr\}.
\end{eqnarray*}

Let $\varphi$ be a function such that
$\varphi (x)\geq I_{(-\infty,0]}(x)$
for all $x\in {\Bbb R}$ and $\varphi $ satisfies the Lipschitz
condition with constant $L({\varphi})$.

\begin{theorem}
Let ${\cal F}:={\rm conv}({\cal H}),$ where ${\cal H}$ is a
class of measurable functions from $(S,{\cal A})$ into ${\Bbb R}.$
For all $t >0,$
$$
{\Bbb P}\Bigl\{\exists f\in {\cal F}:
P\{\tilde f\leq 0\} > \inf_{\delta \in (0,1]}
\Bigl[P_n \varphi ({{\tilde f}\over {\delta}})+
{{8L(\varphi)}\over {\delta}}R_n({\cal H})+
$$
$$
\Bigl(\frac{\log\log_2 (2\delta^{-1})}{n}\Bigr)^{1/2}\Bigr]+
\frac{t}{\sqrt{n}}
\Bigr\}\leq
2\exp \{-2t^2\}.
$$
\end{theorem}

{\bf Proof}.
Since ${\cal F}:={\rm conv}({\cal H}),$ where ${\cal H}$ is a
class of measurable functions from $(S,{\cal A})$ into ${\Bbb R},$
we have
\begin{eqnarray*}
&&
R_n({\cal F})=
{\Bbb E}\|n^{-1}\sum_{i=1}^n \eps_i\delta_{X_i}\|_{\cal F}
\\
&&
={\Bbb E}\sup\Bigl\{|n^{-1}\sum_{i=1}^n \eps_i f_N(X_i)|:
f_N=\sum_{j=1}^{N} w_j^N h_j^N,
w_j^N\geq 0, \sum_{j=1}^{N} w_j^N =1,\ h_j^N\in {\cal H},\ N\geq 1
\Bigr\}
\\
&&
={\Bbb E}\|n^{-1}\sum_{i=1}^n \eps_i\delta_{X_i}\|_{\cal H}=R_n({\cal H}).
\end{eqnarray*}
It follows that $R_n(\tilde {\cal F})=R_n({\cal H}),$ and Theorem 2
implies  the result.
\qed

In the voting methods of combining the classifiers
(such as boosting, bagging (Breiman (1996)), etc.),
a classifier produced at each iteration is a convex combination 
$ f_{\cal S}\in {\rm conv}({\cal H})$
of simple base classifiers from the class ${\cal H}$ 
($f_{\cal S}$ depends on the training sample ${\cal S}:=((X_1,Y_1),\dots, (X_n,Y_n))$). 
The bound of Theorem 12 implies
that for a given $\alpha \in (0,1)$ with probability at least $1-\alpha$
$$
P\{\tilde f_{\cal S}\leq 0\}\leq
\inf_{\delta\in (0,1]}\Bigl[P_n\{\tilde f_{\cal S}\leq \delta\} +
{{8}\over {\delta}}R_n({\cal H})+
\Bigl(\frac{\log\log_2 (2\delta^{-1})}{n}\Bigr)^{1/2}\Bigr]+
\frac{t_{\alpha}}{\sqrt{n}},
$$
where $t_{\alpha}:=\sqrt{{1\over 2}\log{2\over {\alpha}}}.$
In particular, if ${\cal H}$ is a VC--class of classifiers $h:S\mapsto \{-1,1\}$
(which means that the class of sets $\{\{x: h(x)=+1\}: h\in {\cal H}\}$ is
a Vapnik--Chervonenkis class) with VC--dimension $V({\cal H})$, we have with
some constant $C>0$
$$
R_n({\cal H})\leq C\sqrt{{V({\cal H})}\over {n}}.
$$
This implies that with probability at least $1-\alpha$
$$
P\{\tilde f_{\cal S}\leq 0\}\leq \inf_{\delta\in (0,1]}
\Bigl[P_n\{\tilde f_{\cal S} \leq \delta\} +
{C\over {\delta}}\sqrt{{V({\cal H})}\over {n}}+
\Bigl(\frac{\log\log_2(2\delta^{-1})}{n}\Bigr)^{1/2}
\Bigr]
+\sqrt{{1\over {2n}}\log{2\over {\alpha}}},
$$
which slightly improves the main bound
of the paper of Schapire, Freund, Bartlett and Lee (1998),
which has a factor $\log(n/V({\cal H}))$ in front of the term
$C\delta^{-1}(V({\cal H})/n)^{1/2}.$

{\bf Example}. In this example we consider a popular boosting algorithm called
AdaBoost. At the beginning (at the first iteration) AdaBoost assigns uniform
weights $w_j^{(1)}=n^{-1}$ to the labeled observations $(X_1,Y_1),\dots ,(X_n,Y_n).$
At each iteration the algorithm updates the weights. Let $w^{(k)}=(w_1^{(k)},\dots ,w_n^{(k)})$
denote the vector of weights at $k$-th iteration. Let $P_{n,w^{(k)}}$ be the weighted
empirical measure on the $k$-th iteration:
$$
P_{n,w^{(k)}}:=\sum_{i=1}^n w_i^{(k)} \delta_{(X_i,Y_i)}.
$$
AdaBoost calls iteratively a base learning algorithm (called "weak learner")
that returns at $k$-th iteration a classifier $h_k\in {\cal H}$ and
computes the weighted training error of $h_k:$
$$
e_k := P_{n,w^{(k)}}\{y\neq h_k\}.
$$
(In fact, the weak learner attempts to find a classifier with small enough
weighted training error, at least such that $e_k\leq 1/2$).
Then the weights are updated according to the rule
$$
w_j^{(k+1)}:={{w_j^{(k)}\exp\{-Y_j \alpha_k h_k(X_j)\}}\over {Z_k}},
$$
where
$$
Z_k := \sum_{j=1}^n w_j^{(k)}\exp\{-Y_j \alpha_k h_k(X_j)\}
$$
and
$$
\alpha_k :={1\over 2}\log{{1-e_k}\over {e_k}}.
$$
After $N$ iterations AdaBoost outputs a classifier
$$
f_{\cal S} (x):={{\sum_{k=1}^N \alpha_k h_k(x)}\over {\sum_{k=1}^N \alpha_k}}.
$$
The above bounds, of course, apply to this classifier since
$ f_{\cal S}\in {\rm conv}({\cal H}).$ Another way to use Theorem 12
in the case of this example is to choose a decreasing function $\varphi ,$
satisfying all the conditions of Theorem 12 with
$L(\varphi)=1$ and such that
$\varphi (u)\leq e^{-u}$ for all $u\in {\Bbb R}.$ It is easy
to see that such a choice is possible. Let us also set
$$
\delta :={1\over {\sum_1^N \alpha_k}}\bigwedge 1.
$$
Then it is not hard to check that
$$
\varphi \Bigl({{y\sum_1^N \alpha_k h_k(x)}\over {\delta \sum_1^N\alpha_k}}\Bigr)\leq
\varphi (y\sum_1^N \alpha_k h_k(x))\leq
\exp \{-y\sum_1^N \alpha_k h_k(x)\}.
$$
Therefore
$$
P_n \varphi (\frac{\tilde f_{\cal S}}{\delta})\leq
P_n \exp\{-y\sum_1^N \alpha_k h_k(x)\}.
$$
A simple (and well known in the literature on boosting, see e.g.
Schapire, Freund, Bartlett and Lee (1998)) computation
shows that
$$
P_n \exp\{-y\sum_1^N \alpha_k h_k(x)\}=\prod_{k=1}^N 2\sqrt{e_k(1-e_k)}.
$$
We also have
$$
\sum_{k=1}^N \alpha_k = \log{\prod_{k=1}^N \sqrt{{1-e_k}\over {e_k}}}.
$$
It follows now from the bound of Theorem 12 that with
probability at least $1-\alpha$
\begin{eqnarray*}
&&
P\{\tilde f_{\cal S}\leq 0\} \leq  \prod_{k=1}^N 2\sqrt{e_k(1-e_k)} +
8\Bigl(\log{\prod_{k=1}^N \sqrt{{1-e_k}\over {e_k}}}\bigvee 1\Bigr)
R_n({\cal H})
\\
&&
+\biggl(\frac{\log\log_2
\bigl(2\bigr(\log{\prod_{k=1}^N \sqrt{{1-e_k}\over {e_k}}}\bigvee 1\bigl)
\bigr)}{n}\biggr)^{1/2}
+\sqrt{{1\over {2n}}\log{2\over {\alpha}}}.
\end{eqnarray*}

The results of Section 3 provide some improvements of the above
bounds on generalization error of convex combinations of base
classifiers. To be specific, consider the case when ${\cal H}$
is a VC-class of classifiers. Let $V:=V({\cal H})$ be its
VC-dimension.
A well known bound
on the entropy of the convex hull of a VC-class
(see van der Vaart and Wellner (1996), p. 142)
implies that
$$
H_{d_{P_n,2}}({\rm conv}({\cal H});u)\leq
\sup_{Q\in {\cal P}(S)}H_{d_{Q,2}}({\rm conv}({\cal H});u)\leq
D u^{-\frac{2(V-1)}{V}}.
$$
[The bound on the entropy of a convex hull goes back to Dudley;
the precise value of the exponent was given by Ball and Pajor,
van der Vaart and Wellner, Carl; in the case of the convex hull of
a VC-class, the above bound relies also on Haussler's improvement
of Dudley's original bound on the entropy of a VC-class.
See the discussion in the books of van der Vaart and Wellner (1996)
and Dudley (1999) and references therein.]
It immediately follows from Theorem 5 that for all
$\gamma \geq \frac{2(V-1)}{2V-1}$ and for some constants $C, B$
$$
{\Bbb P}\Bigl\{\exists f\in {{\rm conv}({\cal H})}:
P\{\tilde f\leq 0\}>\frac{C}
{n^{1-\gamma/2}\hat \delta_n(\gamma ;f)^{\gamma}}
\Bigr\}\leq
B\log_2\log_2{n}\exp\Bigl\{-\frac{1}{2}n^{\frac{\gamma}{2}}\Bigr\},
$$
where
$$
\hat \delta_n(\gamma ;f):=
\sup\Bigl\{\delta\in (0,1):
\delta^{\gamma }
P_n\{(x,y):y f(x)\leq \delta\}\leq n^{-1+\frac{\gamma}{2}}
\Bigr\}.
$$
This shows that in the case when the VC-dimension of the
base is relatively small the generalization error of boosting
and some other convex combinations of simple classifiers obtained
by various versions of voting methods becomes better than it was
suggested by the bounds of Schapire, Freund, Bartlett and Lee (1998).
One can also conjecture, based on the bounds of Section 3,
that outstanding generalization ability
of these methods observed in numerous experiments can be related
not only to the fact that they produce large margin classifiers,
but also to the fact that the combined classifier belongs to a subset of the
whole convex hull for which the random entropy $H_{d_{P_n,2}}$
is much smaller than for the whole convex hull.

Finally, it is worth mentioning that the bounds in terms
of the so called margin cost functions (see e.g. Mason, Bartlett
and Baxter (1999), Mason, Baxter, Bartlett and Frean (1999))
easily follow from Theorem 1. 
Namely, Theorem 1 implies that with probability at least
$1-\alpha$
$$
P\{\tilde f_{\cal S}\leq 0\}\leq \inf_{N\geq 
1}\Bigl[P_n\varphi_N(\tilde f_{\cal S}) +
CL_N \sqrt{{V({\cal H})}\over {n}}+
\Bigl(\frac{\log N}{n}\Bigr)^{1/2}
\Bigr]
+\sqrt{{1\over {2n}}\log{2\over {\alpha}}},
$$
where $\{\varphi_N\}$ is  any  sequence of  Lipschitz  cost
functions such that $\varphi_N (x)\geq I_{(-\infty,0]}(x)$ for all
$x\in {\Bbb R}, N\geq 1$ and $L_N$ is a Lipschitz constant of
$\varphi_N.$

\section{Bounding the generalization error in neural network learning}

We turn now to the applications of the bounds of previous
section in neural network learning. We start with the description
of the class of feedforward neural networks for which the bounds
on the generalization error will be proved.
Let ${\cal H}$ be a class of measurable functions from
$(S,{\cal A})$ into ${\Bbb R}$ (base functions).
Consider an acyclic directed graph $G.$ Suppose that
$G$ has a unique vertex $v_i$ (input) that has no incoming
edges and a unique vertex $v_o$ (output) that has one outcoming
edge. The vertices (nodes) of the graph will be called \it neurons. \rm
Suppose the set $V$ of all the neurons
is divided into layers
$$
V=\{v_i\}\cup \bigcup_{j=0}^l V_j,
$$
where $l\geq 0$ and $V_l=\{v_o\}.$
The neurons $v_i, v_o$ are called the input and the output neurons,
respectively. The neurons of the layer $V_0$ will be called the
base neurons.
Suppose also that the inputs of the base neurons
are the outputs of the input neuron.
Suppose also that the inputs
of the neurons of the layer $V_j, j\geq 1$
are the ouputs of the neurons from the set
$
\bigcup_{k=0}^{j-1} V_k.
$
To define the network, we will assign the labels to the
neurons the following way. Each of the base neurons is labeled
by a function from the base class ${\cal H}.$ Each neuron
of the $j$th layer $V_j,$ where $j\geq 1,$ is labeled by
a vector $w:=(w_1,\dots ,w_n)\in {\Bbb R}^n,$ where $n$
is the number of inputs of the neuron. $w$ will be called
the vector of weights of the neuron.
 
Given a Borel function $\sigma $ from ${\Bbb R}$ into $[-1,1]$
(a sigmoid) and a vector $w:=(w_1,\dots ,w_n)\in {\Bbb R}^n,$
let
$$
N_{\sigma ,w}:{\Bbb R^n}\mapsto {\Bbb R},\ N_{\sigma ,w}(u_1,\dots ,u_n):=\sigma (\sum_{i=1}^n w_j u_j).
$$
For $w\in {\Bbb R}^n,$
$$
\|w\|_{{\ell}_1} := \sum_{i=1}^n |w_i|.
$$
Let $\sigma_j: j\geq 1$ be functions from ${\Bbb R}$ into $[-1,1],$ satisfying the Lipschitz conditions:
$$
|\sigma_j (u)-\sigma_j (v)|\leq L_j |u-v|,\ u,v\in {\Bbb R}.
$$

The network works the following way. The input neuron inputs
an instance $x\in S.$ A base neuron computes the value
of the base function (it is labeled with) on this instance
and outputs the value through its output edges.
A neuron in $j$th layer ($j\geq 1$) computes and outputs
through its output edges the value $N_{\sigma_j,w}(u_1,...,u_n)$
(where $u_1,\dots ,u_n$ are the values of the inputs of the neuron).
The network outputs the value $f(x)$ (of a function $f$ it computes)
through the output edge.

We denote ${\cal N}_l$ the set of all such networks.
We call ${\cal N}_l$ the class of feedforward neural networks with base ${\cal H}$
and $l$ layers of neurons (and with sigmoids $\{\sigma_j\}$).
Let
${\cal N}_{\infty}:=\bigcup_{j=0}^{\infty } {\cal N}_j.$
Define ${\cal H}_0 := {\cal H},$ and then recursively
$$
{\cal H}_j := \Bigl\{N_{\sigma_j, w}(h_1,\dots ,h_n): n\geq 0, h_i\in {\cal H}_{j-1},\ w\in {\Bbb R}^n\Bigr\}
\bigcup {\cal H}_{j-1}.
$$
Denote
$
{\cal H}_{\infty}:=\bigcup_{j=0}^{\infty } {\cal H}_j.
$
Clearly, ${\cal H}_{\infty}$ includes all the functions
computable by feedforward neural networks with base ${\cal H}.$

Let $\{A_j\}$ be a sequence of positive numbers.
We also define recursively classes of functions
computable by feedforward neural networks with restrictions on the
weights of neurons:
$$
{\cal H}_j (A_1,\dots ,A_j):=
$$
$$
:=\Bigl\{N_{\sigma_j, w}(h_1,\dots ,h_n): n\geq 0, h_i\in {\cal H}_{j-1}(A_1,\dots ,A_{j-1}),\
w\in {\Bbb R}^n, \|w\|_{{\ell}_1}\leq A_j \Bigr\}\bigcup
$$
$$
\bigcup {\cal H}_{j-1}(A_1,\dots ,A_{j-1}).
$$
Clearly,
$$
{\cal H}_j := \bigcup \Bigl\{{\cal H}_j(A_1,\dots ,A_j): A_1,\dots ,A_j<+\infty \Bigr\}.
$$
As in the previous section, let $\varphi$ be a function such that
$\varphi (x)\geq I_{(-\infty,0]}(x)$
for all $x\in {\Bbb R}$ and $\varphi $ satisfies the Lipschitz
condition with constant $L({\varphi})$.

We start with the following result.

\begin{theorem}
For all $t >0$ and for all $l\geq 1$
$$
{\Bbb P}\Bigl\{\exists f\in {\cal H}_l(A_1,\dots ,A_l):
P\{\tilde f\leq 0\} > \inf_{\delta \in (0,1]}
\Bigl[P_n \varphi ({{\tilde f}\over {\delta}})+
{{2\sqrt{2\pi}L(\varphi)}\over {\delta}} \prod_{j=1}^l (2L_j A_j+1)
G_n({\cal H})\Bigr]
$$
$$
+\frac{t+2}{\sqrt{n}}
\Bigr\}\leq
2\exp \{-2t^2\}.
$$
\end{theorem}

{\bf Proof}. We apply Theorem 2 to the class ${\cal F}={\cal H}_l(A_1,\dots ,A_l)=:{\cal H}_l^{\prime},$
which gives for all $t>0$
$$
{\Bbb P}\Bigl\{\exists f\in {\cal H}_l^{\prime}:
P\{\tilde f\leq 0\} >
\inf_{\delta \in [0,1]}\Bigl[P_n \varphi ({{\tilde f}\over {\delta}})+
{{2\sqrt{2\pi}L(\varphi)}\over {\delta}}G_n({\cal H}_l^{\prime})+
\Bigl(\frac{\log\log_2 (2\delta^{-1})}{n}\Bigr)^{1/2}\Bigr]+
\frac{t+2}{\sqrt{n}}
\Bigr\}
$$
$$
\leq 2\exp \{-2t^2\}.
$$
Thus, it's enough to show that
$$
G_n({\cal H}_l^{\prime})={\Bbb E}\|n^{-1}\sum_{i=1}^n g_i\delta_{X_i}\|_{{\cal H}_l^{\prime}}\leq
\prod_{j=1}^l (2L_j A_j+1)
{\Bbb E}\|n^{-1}\sum_{i=1}^n g_i\delta_{X_i}\|_{\cal H}.
$$
To this end, note that
\begin{equation}
{\Bbb E}\|n^{-1}\sum_{i=1}^n g_i\delta_{X_i}\|_{{\cal H}_l^{\prime}}\leq
{\Bbb E}\|n^{-1}\sum_{i=1}^n g_i\delta_{X_i}\|_{{\cal G}_{l}}+
{\Bbb E}\|n^{-1}\sum_{i=1}^n g_i\delta_{X_i}\|_{{\cal H}_{l-1}^{\prime}},
\label{e6.1}
\end{equation}
where
$$
{\cal G}_l :=
\Bigl\{N_{\sigma_l, w}(h_1,\dots ,h_n): n\geq 0, h_i\in {\cal H}_{l-1}(A_1,\dots ,A_{l-1}),\
w\in {\Bbb R}^n, \|w\|_{{\ell}_1}\leq A_l \Bigr\}.
$$
Consider two Gaussian processes
$$
Z_1(f):= n^{-1/2}\sum_{i=1}^n g_i (\sigma_l \circ f)(X_i)
$$
and
$$
Z_2(f):= L_l n^{-1/2}\sum_{i=1}^n g_i f(X_i),
$$
where
$$
f\in \Bigl\{\sum_{i=1}^n w_i h_i: n\geq 0, h_i\in {\cal H}_{l-1}^{\prime},\ w\in {\Bbb R}^n,\ \|w\|_{{\ell}_1}\leq A_l
\Bigr\}=: {\cal G}_l^{\prime}.
$$
We have
\begin{eqnarray*}
&&
{\Bbb E}_g|Z_1(f)-Z_1(h)|^2 =
n^{-1}\sum_{i=1}^n |\sigma_l(f(X_i)-\sigma_l(h(X_i))|^2
\\
&&
\leq
L_l^2 n^{-1}\sum_{i=1}^n |f(X_i)-h(X_i)|^2=
{\Bbb E}_g|Z_2(f)-Z_2(h)|^2.
\end{eqnarray*}
By Slepian's Lemma (see Ledoux and Talagrand (1991)), we get
\begin{eqnarray}
{\Bbb E}_g\|n^{-1}\sum_{i=1}^n g_i\delta_{X_i}\|_{{\cal G}_{l}}= 
n^{-1/2}{\Bbb E}_g\|Z_1\|_{{\cal G}_l^{\prime}}
\leq
2 n^{-1/2}{\Bbb E}_g\|Z_2\|_{{\cal G}_l^{\prime}}=
2 L_l
{\Bbb E}_g\|n^{-1}\sum_{i=1}^n g_i\delta_{X_i}\|_{{\cal G}_{l}^{\prime }}.
\label{e6.2}
\end{eqnarray}
Since ${\cal G}_l^{\prime }=A_l {\rm conv}_s({\cal H}_{l-1})$
[here ${\rm conv}_s({\cal G})$ denotes closed symmetric convex hull
of a class ${\cal G}$, i.e. closed convex hull of the class
${\cal H}\cup -{\cal H}$],
it is easy to get that
\begin{equation}
{\Bbb E}\|n^{-1}\sum_{i=1}^n g_i\delta_{X_i}\|_{{\cal G}_{l}^{\prime }}=
A_l
{\Bbb E}\|n^{-1}\sum_{i=1}^n g_i\delta_{X_i}\|_{{\cal H}_{l-1}}.
\label{e6.3}
\end{equation}
It follows from the bounds (\ref{e6.1})--(\ref{e6.3}) that
$$
{\Bbb E}\|n^{-1}\sum_{i=1}^n g_i\delta_{X_i}\|_{{\cal H}_{l}}\leq
(2L_l A_l+1)
{\Bbb E}\|n^{-1}\sum_{i=1}^n g_i\delta_{X_i}\|_{{\cal H}_{l-1}}.
$$
The result now follows by induction.
\qed

{\bf Remark}. It can be shown that in the case of multilayer
perceptrons (in which the neurons in each layer are linked only
to the neurons in the previous layer) the factor
$\prod_{j=1}^l (2L_j A_j+1)$ in the bound of the theorem can
be replaced by $\prod_{j=1}^l (2L_j A_j).$ If the sigmoids
are odd functions, the same factor in the case of general feedforward
architecture of the network becomes $\prod_{j=1}^l (L_j A_j+1),$
and in the case of multilayer perceptrons
$\prod_{j=1}^l L_j A_j.$
Bartlett (1998) obtained a bound similar to the first inequality
of Theorem 13 for a more special
class ${\cal H}$ and with larger constants. In the case when $A_j \equiv A, L_j \equiv L$ (the case
considered by Bartlett) the expression in the right hand side of his bound includes ${{(AL)^{l(l+1)/2}}\over {\delta^l}},$
which is replaced in our bound by ${{(AL)^l}\over {\delta}}.$ These improvement can be substantial in applications,
since the above quantities play the role of complexity penalties.

Given a neural network $f\in {\cal N}_{\infty},$ let
$$
{\ell}(f):=\min \{j\geq 1: f\in {\cal N}_j\}.
$$
Let $\{b_k\}$ be a sequence of nonnegative numbers.
For a number $k, 1\leq k\leq {\ell}(f),$
let ${V}_k(f)$ denote the
set of all neurons
of layer $k$ in the graph representing $f.$ Denote
$$
W_k(f) := \max_{N\in {V}_k(f)} \|w^{(N)}\|_{{\ell}_1}\bigvee b_k,\ k=1, 2, \dots , \ell (f),
$$
and let
$$
\Lambda (f):= \prod_{k=1}^{\ell (f)} (4 L_k W_k(f)+1),
$$
$$
\Gamma_{\alpha} (f):= \sum_{k=1}^{\ell (f)} \sqrt{{{\alpha}\over 2}\log (2+|\log_2 W_k(f)|)},
$$
where $\alpha >0$ is a number such that $\zeta (\alpha)<3/2,$
$\zeta$ being the Riemann zeta-function:
$$
\zeta (\alpha) := \sum_{k=1}^{\infty }k^{-\alpha}.
$$

\begin{theorem}
For all $t >0$ and for all
$\alpha >0$ such that $\zeta (\alpha )<3/2,$ the following bounds
hold:
\begin{eqnarray*}
&&
{\Bbb P}\Bigl\{\exists f\in {\cal H}_{\infty}:
P\{\tilde f\leq 0\} > \inf_{\delta \in (0,1)}\Bigl[P_n \varphi ({{\tilde f}\over {\delta}})+
+{{2\sqrt{2\pi}L(\varphi)}\over {\delta}} \Lambda (f) G_n({\cal H})
\\
&&
+\Bigl(\frac{\log\log_2 (2\delta^{-1})}{n}\Bigr)^{1/2}
\Bigr]+
{{\Gamma_{\alpha} (f)+t+2}\over {\sqrt {n}}}
\Bigr\}\leq 2(3-2\zeta (\alpha))^{-1}\exp \{-2t^2\}.
\end{eqnarray*}
\end{theorem}

{\bf Proof}. With a little abuse of notations, we write $f$
for both the neural network and the function it computes.
Denote
$$
\Delta_k := \cases {
                 [2^{k-1},2^k) & for $k\in {\Bbb Z}, k\neq 0, 1$\cr
                 [1/2,2) & for $k=1.$\cr
               }
$$
The conditions $\ell (f)=l$ and
$$
W_j(f)\in {\Delta_{k_j}},\ k_j\in {\Bbb Z}\setminus \{0\},\ j=1,\dots ,l
$$
easily imply that
$$
\Lambda (f)\geq \prod_{j=1}^l (2L_j 2^{k_j}+1),\,\,\,
\Gamma_{\alpha}(f)\geq
\sum_{j=1}^l \sqrt {{{\alpha }\over 2}\log (|k_j|+1)} 
$$
and also that
$f\in {\cal H}_l (2^{k_1},\dots ,2^{k_l}).$
Therefore, the following bounds hold:
$$
{\Bbb P}\Bigl\{\exists f\in {\cal H}_{\infty}:
P\{\tilde f\leq 0\} > \inf_{\delta \in (0,1)}
\Bigl[P_n \varphi ({{\tilde f}\over {\delta}})+
{{2\sqrt{2\pi}L(\varphi)}\over {\delta}} \Lambda (f) G_n({\cal H})
$$
$$
+\Bigl(\frac{\log\log_2 (2\delta^{-1})}{n}\Bigr)^{1/2}
\Bigr]+
{{\Gamma_{\alpha} (f)+t+2}\over {\sqrt {n}}}
\Bigr\}
$$
$$
\leq
\sum_{l=0}^{\infty}\sum_{k_1\in {\Bbbb Z}\setminus \{0\}} \dots \sum_{k_l\in {\Bbbb Z}\setminus \{0\}}
{\Bbb P}\Bigl\{\exists f\in {\cal H}_{\infty}\bigcap
\bigl\{f:\ell (f)=l,\ W_j(f)\in \Delta_{k_j},\ j=1,\dots ,l\bigr\}:
$$
$$
P\{\tilde f\leq 0\} > \inf_{\delta \in (0,1)}
\Bigl[P_n \varphi ({{\tilde f}\over {\delta}})
+{{2\sqrt{2\pi}L(\varphi)}\over {\delta}} \Lambda (f) G_n({\cal H})
$$
$$
+\Bigl(\frac{\log\log_2 (2\delta^{-1})}{n}\Bigr)^{1/2}
\Bigr]+
{{\Gamma_{\alpha} (f)+t+2}\over {\sqrt {n}}}
\Bigr\}
$$
$$
\leq
\sum_{l=0}^{\infty}\sum_{k_1\in {\Bbbb Z}\setminus \{0\}} \dots \sum_{k_l\in {\Bbbb Z}\setminus \{0\}}
{\Bbb P}\Bigl\{\exists f\in {\cal H}_l(2^{k_1},\dots, 2^{k_l}):
P\{\tilde f\leq 0\} > \inf_{\delta \in (0,1)}\Bigl[P_n \varphi 
({{\tilde f}\over {\delta}})
$$
$$
+{{2\sqrt{2\pi}L(\varphi)}\over {\delta}}
\prod_{j=1}^l(2L_j 2^{k_j}+1) G_n({\cal H})
+\Bigl(\frac{\log\log_2 (2\delta^{-1})}{n}\Bigr)^{1/2}
\Bigr]
$$
$$
+{{\sum_{j=1}^l \sqrt{{{\alpha }\over 2}\log (|k_j|+1)} + t+2}\over {\sqrt {n}}}
\Bigr\}.
$$
Using the bound of Theorem 13, we obtain
$$
{\Bbb P}\Bigl\{\exists f\in {\cal H}_{\infty}:
P\{\tilde f\leq 0\} > \inf_{\delta \in (0,1)}\Bigl[P_n \varphi ({{\tilde f}\over {\delta}})+
{{2\sqrt{2\pi}L(\varphi)}\over {\delta}} \Lambda (f) G_n({\cal H})
$$
$$
+\Bigl(\frac{\log\log_2 (2\delta^{-1})}{n}\Bigr)^{1/2}
\Bigr]+
{{\Gamma_{\alpha} (f)+t+2}\over {\sqrt {n}}}
\Bigr\}
$$
$$
\leq
\sum_{l=0}^{\infty}\sum_{k_1\in {\Bbbb Z}\setminus \{0\}} \dots \sum_{k_l\in {\Bbbb Z}\setminus \{0\}}
2\exp \{-2(\sum_{j=1}^l \sqrt{{{\alpha }\over 2}\log (|k_j|+1)} + t)^2\}
$$
$$
\leq
\sum_{l=0}^{\infty}\sum_{k_1\in {\Bbbb Z}\setminus \{0\}} \dots \sum_{k_l\in {\Bbbb Z}\setminus \{0\}}
2\exp \{-\sum_{j=1}^l \alpha\log (|k_j|+1)-2t^2\}
$$
$$
=2\sum_{l=0}^{\infty}\sum_{k_1\in {\Bbbb Z}\setminus \{0\}} \dots \sum_{k_l\in {\Bbbb Z}\setminus \{0\}}
\prod_{j=1}^l (|k_j|+1)^{-\alpha}\exp \{-2t^2\}
$$
$$
=2\sum_{l=0}^{\infty }\prod_{j=1}^l (2\sum_{k=2}^{\infty} k^{-\alpha}) \exp\{-2t^2\}=
2\sum_{l=0}^{\infty}[2(\zeta (\alpha)-1)]^l \exp\{-2t^2\}
$$
$$
=2(3-2\zeta (\alpha))^{-1}\exp \{-2t^2\}
$$
which yields the bound of the theorem.
\qed

It follows, in particular, that for any classifier 
$ f_{\cal S}\in {\cal H}_{\infty},$
based on the training data 
${\cal S}:=((X_1,Y_1),\dots ,(X_n,Y_n)),$ we have
$$
{\Bbb P}\Bigl\{
P\{\tilde f_{\cal S}\leq 0\} > \inf_{\delta \in (0,1)}
\Bigl[P_n \varphi ({{\tilde f_{\cal S}}\over {\delta}})+
{{2\sqrt{2\pi}L(\varphi)}\over {\delta}} \Lambda (f_{\cal S}) G_n({\cal H})+
$$
$$
+\Bigl(\frac{\log\log_2 (2\delta^{-1})}{n}\Bigr)^{1/2}
\Bigr]+
{{\Gamma_{\alpha} (f_{\cal S})+t+2}\over {\sqrt {n}}}
\Bigr\}\leq 2(3-2\zeta (\alpha))^{-1}\exp \{-2t^2\}.
$$

Next we consider a method of complexity penalization in neural network learning
based on the penalties that depend on $\ell_1$-norms of the vectors of weights
of the neurons.
Suppose that $f_{\cal S}$ is the neural network from 
${\cal F}\subset {\cal H}_{\infty}$ that minimizes
the penalized training error
$$
f_{\cal S} :=
{\rm arg min}_{f\in {\cal F}}
\inf_{\delta \in (0,1]}\Bigl[P_n (\{\tilde f\leq \delta \})+
{{2\sqrt{2\pi}}\over {\delta}} \Lambda (f) G_n({\cal H})
+\Bigl(\frac{\log\log_2 (2\delta^{-1})}{n}\Bigr)^{1/2}
\Bigr]+
{{\Gamma_{\alpha} (f)}\over {\sqrt{n}}}
$$
$$
={\rm arg min}_{f\in {\cal F}} \Bigl[P_n(\{\tilde f\leq 0\})+
\inf_{\delta \in (0,1]} \hat \pi_n (f;\delta)\Bigr],
$$
where the quantity
$\inf_{\delta \in (0,1]}\hat \pi_n(f;\delta)$
plays the role of the complexity penalty,
$$
\hat \pi_n (f;\delta):=P_n (\{0<\tilde f\leq \delta \})+
\Psi_n (f;\delta),
$$
$$
\Psi_n (f;\delta):=
{{2\sqrt{2\pi}}\over {\delta}} \Lambda (f) G_n({\cal H})
+\Bigl(\frac{\log\log_2 (2\delta^{-1})}{n}\Bigr)^{1/2}+
{{\Gamma_{\alpha} (f)}\over {\sqrt{n}}}.
$$
We define a distribution dependent version of this data dependent penalty
as $\inf_{\delta \in (0,1]}\pi_n(f;\delta),$
where
$$
\pi_n (f;\delta):=P(\{0<\tilde f\leq 2\delta \})+
2\Psi_n(f;\delta).
$$
The first inequality of the next theorem provides
an upper confidence bound on the generalization
error of the classifier $f_{\cal S}.$ The second bound
is an "oracle inequality" that shows that the estimate $f_{\cal S}$
obtained by the above method possess some optimality property (see
Johnstone (1998), Barron, Birg\' e and
Massart (1999) for a general approach to penalization and oracle inequalities
in nonparametric statistics).

\begin{theorem}
For all $t >0$ and for all
$\alpha >0$ with $\zeta (\alpha)<3/2,$
the following bounds hold:
$$
{\Bbb P}\Bigl\{
P\{{\tilde {f_{\cal S}}}\leq 0\} >
\inf_{f\in {\cal F}} \bigl[P_n\{\tilde f\leq 0\}+
\inf_{\delta \in (0,1]}\hat \pi_n(f;\delta)
\bigr]+{{t+2}\over {\sqrt{n}}}
\Bigr\}\leq 2(3-2\zeta (\alpha))^{-1} \exp \{-2t^2\}
$$
and
\begin{eqnarray*}
{\Bbb P}\Bigl\{
P\{\tilde {f_{\cal S}}\leq 0\}-\inf_{g\in {\cal F}}P\{\tilde g\leq 0\} 
&>& 
\inf_{f\in {\cal F}}
\Bigl[P\{\tilde f\leq 0\}-\inf_{g\in {\cal F}}P\{\tilde g\leq 0\}
+\inf_{\delta \in (0,1]}\pi_n(f;\delta)
\Bigr]
\\
&+&
{{2t+4}\over {\sqrt{n}}}\Bigr\}
\leq 4(3-2\zeta(\alpha))^{-1}\exp \{-2t^2\}.
\end{eqnarray*}
\end{theorem}

{\bf Proof}. The first bound follows from Theorem 14 and the definition of
the estimate $\tilde f_{\cal S}.$
To prove the second bound, we repeat the proof of Theorems 1, 2 to show that
for any class ${\cal F}^{\prime}$
$$
{\Bbb P}\Bigl\{\exists f\in {\cal F}^{\prime}\ \exists \delta \in (0,1]:
P_n\{\tilde f\leq \delta\} >
\Bigl[P\varphi ({{\tilde f-\delta}\over {\delta}})
+{{2\sqrt{2\pi}}\over {\delta}} G_n({\cal F}^{\prime})+
$$
$$
+\Bigl(\frac{\log\log_2 (2\delta^{-1})}{n}\Bigr)^{1/2}
\Bigr]+
{{t+2}\over {\sqrt {n}}}
\Bigr\}\leq 2\exp \{-2t^2\}.
$$
The argument that led to Theorems 13 and 14 shows that
$$
{\Bbb P}\Bigl\{\exists f\in {\cal F}\ \exists \delta \in (0,1]:
P_n\{\tilde f\leq \delta\} >
\Bigl[P\{\tilde f\leq 2\delta\}
+{{2\sqrt{2\pi}}\over {\delta}} \Lambda (f) G_n({\cal H})
$$
$$
+\Bigl(\frac{\log\log_2 (2\delta^{-1})}{n}\Bigr)^{1/2}
+{{\Gamma_{\alpha} (f)}\over {\sqrt {n}}} +
{{t+2}\over {\sqrt {n}}}\Bigr]
\Bigr\}
\leq 2(3-2\zeta(\alpha))^{-1}\exp \{-2t^2\}.
$$
If now
$$
\inf_{f\in {\cal F}}\inf_{\delta \in (0,1]}\Bigl[P_n (\{\tilde f\leq \delta\})+
\Psi_n(f;\delta)\Bigr]+{{t+2}\over {\sqrt{n}}}
$$
$$
>\inf_{f\in {\cal F}} \inf_{\delta \in (0,1]}
\Bigl[P\{\tilde f\leq 2\delta\}+
2\Psi_n(f;\delta)\Bigr]+{{2t+4}\over {\sqrt{n}}},
$$
then
$$
\exists f\in {\cal F}\ \exists \delta \in (0,1]:
P_n\{\tilde f\leq \delta\} >
\Bigl[P\{\tilde f\leq 2\delta\}+
\Psi_n(f;\delta)\Bigr]+{{t+2}\over {\sqrt{n}}}.
$$
Combining this with the first bound gives
$$
{\Bbb P}\Bigl\{
P\{\tilde {f_{\cal S}}\leq 0\} >
\inf_{f\in {\cal F}} \inf_{\delta \in (0,1)}
\Bigl[P\{\tilde f\leq 2\delta\}+
2\Psi_n(f;\delta)\Bigr]+{{2t+4}\over {\sqrt{n}}}
\Bigr\}
$$
$$
\leq 4(3-2\zeta(\alpha))^{-1}\exp \{-2t^2\},
$$
which implies the result.
\qed

\vskip 3mm

\section*{Acknowledgments.} The authors are very thankful to Jon Wellner
for reading the manuscript and making a number of comments and suggestions
that improved the paper. We want to thank Evarist Gin\'e
for pointing out the current formulation of Theorem 8.
We also want to thank Fernando Lozano for providing the results of the
experiment presented in Figure 1.

\vskip 1mm

\hfill\break
Department of Mathematics and Statistics\hfill\break
The University of New Mexico\hfill\break
Albuquerque, NM 87131--1141\hfill\break
e-mail: \{vlad, panchenk\}@math.unm.edu;\hfill\break
URL: http:$\backslash\backslash$www.math.unm.edu$\backslash$
\~{}\{vlad, panchenk\}

\end{document}